\documentclass[a4paper,12pt]{article}
\usepackage{amsfonts}
\usepackage[xdvi]{graphicx}
\usepackage[dvips]{color}
\usepackage{amsmath,amssymb}

\newcommand{\R}{{\mathbb R}}
\newcommand{\C}{{\mathbb C}}

\newcommand{\Z}{{\mathbb Z}}

\makeatletter

\@addtoreset{equation}{section}

\newcounter{def}[section]
\renewcommand{\thedef}{\stepcounter{def}\thesection.\@arabic\c@def }

\begin{document}
\setlength{\baselineskip}{24pt}
\begin{center}
\textbf{\LARGE{Weights, Kovalevskaya exponents and the Painlev\'{e} property}}
\end{center}

\setlength{\baselineskip}{14pt}

\begin{center}
Advanced Institute for Materials Research, Tohoku University,\\
Sendai, 980-8577, Japan

\large{Hayato CHIBA} \footnote{E mail address : hchiba@tohoku.ac.jp}
\end{center}
\begin{center}
Oct 13, 2020
\end{center}

\begin{center}
\textbf{Abstract}
\end{center}

Weighted degrees of quasihomogeneous Hamiltonian functions of the Painlev\'{e} equations are investigated.
A tuple of positive integers, called a regular weight, satisfying certain conditions related to
singularity theory is classified.
Each polynomial Painlev\'{e} equation has a regular weight.
Conversely, for $2$ and $4$-dim cases, it is shown that there exists a differential equation
satisfying the Painlev\'{e} property associated with each regular weight.
Kovalevskaya exponents of quasihomogeneous Hamiltonian systems are also investigated by means of
regular weights, singularity theory and dynamical systems theory.
It is shown that there is a one-to-one correspondence between Laurent series solutions
and stable manifolds of the associated vector field obtained by the blow-up of the system.
For $4$-dim autonomous Painlev\'{e} equations,  
the level surface of Hamiltonian functions can be 
decomposed into a disjoint union of stable manifolds.
\\

\noindent \textbf{Keywords}: Painlev\'{e} equation; quasihomogeneous vector field;
regular weight; Kovalevskaya exponent

% \tableofcontents

%%%%%%%%%%%%%%%%%%%%%%%%%%%%%%%%%%%%%%%%%%%%%%%%%%%%%%%%%%%%%%%%%%%%%%%%%%%%%%%%%%%%%%%%%%%%%%%%%%
%%%%%%%%%%%%%%%%%%%%%%%%%%%%%%%%%%%%%%%%%%%%%%%%%%%%%%%%%%%%%%%%%%%%%%%%%%%%%%%%%%%%%%%%%%%%%%%%%%

\section{Introduction}

A differential equation defined on a complex region is said to have the Painlev\'{e} property
if any movable singularity 
(a singularity of a solution which depends on an initial condition) of any solution is a pole.
Painlev\'{e} and his group classified second order ODEs having the Painlev\'{e} property
and found new six differential equations called the Painlev\'{e} equations
$\text{P}_\text{I}, \cdots , \text{P}_\text{VI}$.
Nowadays, it is known that they are written in Hamiltonian forms
\begin{eqnarray*}
(\text{P}_\text{J}): \frac{dq}{dz} = \frac{\partial H_J}{\partial p}, \quad
\frac{dp}{dz} = -\frac{\partial H_J}{\partial q}, \quad J = \text{I}, \cdots, \text{VI}.
\end{eqnarray*}
Among six Painlev\'{e} equations, the Hamiltonian functions of the first, second and fourth
Painlev\'{e} equations are polynomials in both of the independent variable $z$
and the dependent variables $(q,p)$.
They are given by
\begin{eqnarray}
H_{\text{I}} &=& \frac{1}{2}p^2 - 2q^3 - zq, \\
H_\text{II} &=& \frac{1}{2}p^2 - \frac{1}{2}q^4 - \frac{1}{2}zq^2 - \alpha q, \nonumber \\
H_\text{IV} &=& -pq^2 + p^2q - 2pqz - \alpha p + \beta q, \nonumber
\label{1-1}
\end{eqnarray}
respectively, where $\alpha, \beta \in \C$ are arbitrary parameters.

In general, a polynomial $H(x_1, \cdots ,x_n)$ is called a quasihomogeneous polynomial if there 
are positive integers $a_1, \cdots ,a_n$ and $h$ such that
\begin{eqnarray}
H(\lambda ^{a_1}x_1 , \cdots ,\lambda ^{a_n}x_n) = \lambda ^h H(x_1, \cdots ,x_n)
\label{quasi}
\end{eqnarray}
for any $\lambda \in \C$.
A polynomial $H$ is called a semi-quasihomogeneous if $H$ is decomposed into two polynomials as 
$H = H^P + H^N$, where $H^P$ satisfies (\ref{quasi}) and $H^N$ satisfies
\begin{eqnarray*}
H^N(\lambda ^{a_1}x_1 , \cdots ,\lambda ^{a_n}x_n) \sim o(\lambda ^h), \quad |\lambda | \to \infty.
\end{eqnarray*}
The integer $\mathrm{deg}(H) := h$ is called the weighted degree of $H$ with
respect to the weight $\mathrm{deg}(x_1, \cdots ,x_n) := (a_1 ,\cdots ,a_n)$.
$H^P$ and $H^N$ are called the principle part and the non-principle part of $H$, respectively.
The weight of $H$ is determined by the Newton diagram.
Plot all exponents $(r_1, \cdots ,r_n)$ of monomials $x_1^{r_1}x_2^{r_2} \cdots x_n^{r_n}$ included in $H^P$
on the integer lattice in $\R^n$.
If they lie on a unique hyperplane $a_1x_1 + \cdots  + a_nx_n = h$,
then $\mathrm{deg}(H^P) = h$ with respect to the weight $(a_1, \cdots ,a_n)$.
Exponents of monomials included in $H^N$ should be on the lower side of the hyperplane.

The above Hamiltonian functions for $\text{P}_\text{I}, \text{P}_\text{II}$
and $\text{P}_\text{IV}$ are semi-quasihomogeneous.
If we define degrees of variables by $\mathrm{deg} (q,p,z) = (2,3,4)$ for $H_{\text{I}}$,
$\mathrm{deg} (q,p,z) = (1,2,2)$ for $H_{\text{II}}$ and $\mathrm{deg} (q,p,z) = (1,1,1)$ for $H_{\text{IV}}$,
then Hamiltonian functions have the weighted degrees $6,4$ and $3$, respectively, (Table 1)
with $H^N_\text{I} = 0, H^N_\text{II} = - \alpha q$
and $H^N_\text{IV} =  - \alpha p + \beta q$.

The Hamiltonian functions of the third, fifth and sixth Painlev\'{e} equations
are not polynomials in $z$, and their weights include nonpositive integers (Table 1).
They are not treated in this paper, while the analysis of them using weighted projective spaces
is given in \cite{Chi3}.

Higher dimensional Painlev\'{e} equations have not been classified yet, however,
a lot of such equations have been reported in the literature.
A list of four dimensional Painlev\'{e} equations derived from the 
monodromy preserving deformation is given in \cite{KNS, K}.
Lie-algebraic approach is often employed to find new Painlev\'{e} 
equations \cite{DS, FS, GHM, NY}.
Several Painlev\'{e} hierarchies, which are hierarchies of $2n$-dimensional Painlev\'{e} equations,
are obtained by the similarity reductions of soliton equations such as the KdV equation.
Among them, it is known that Hamiltonian functions of the 
the first Painlev\'{e} hierarchy $(\text{P}_\text{I})_n$\cite{Koi, Kud, Shi},
the second-first Painlev\'{e} hierarchy $(\text{P}_\text{II-1})_n$\cite{CJM, CJP, Koi, Kud},  
the second-second Painlev\'{e} hierarchy $(\text{P}_\text{II-2})_n$ and 
the fourth Painlev\'{e} hierarchy $(\text{P}_\text{IV})_n$\cite{GJP, Koi}
can be expressed as polynomials with respect to both of the dependent variables and 
the independent variables.
They are Hamiltonian PDEs of the form
\begin{eqnarray}
\left\{ \begin{array}{l}
\displaystyle \frac{\partial q_j}{\partial z_i} = \frac{\partial H_i}{\partial p_j},\,\,
\frac{\partial p_j}{\partial z_i} = -\frac{\partial H_i}{\partial q_j},\quad
j=1,\cdots ,n;\,\, i=1,\cdots ,n   \\[0.4cm]
H_i = H_i(q_1, \cdots ,q_n, p_1, \cdots ,p_n, z_1, \cdots ,z_n)  \\
\end{array} \right.
\label{1-2}
\end{eqnarray}
consisting of $n$ Hamiltonians $H_1, \cdots ,H_n$ with $n$ independent variables $z_1, \cdots ,z_n$.
When $n=1$, $(\text{P}_\text{I})_1$ and $(\text{P}_\text{IV})_1$ are reduced to the first and fourth
Painlev\'{e} equations, respectively.
Both of $(\text{P}_\text{II-1})_1$ and $(\text{P}_\text{II-2})_1$ coincide with the second Painlev\'{e}
equation, while they are different systems for $n\geq 2$.
When $n=2$, Hamiltonians of $(\text{P}_\text{I})_2$, $(\text{P}_\text{II-1})_2$, 
$(\text{P}_\text{II-2})_2$ and $(\text{P}_\text{IV})_2$ are given by
\begin{equation}
(\text{P}_\text{I})_2 \left\{ \begin{array}{l}
\displaystyle H^{9/2}_1 
=2p_2p_1 + 3p_2^2q_1 + q_1^4 - q_1^2q_2 - q_2^2 - z_1q_1 + z_2(q_1^2 - q_2),  \\
\displaystyle H^{9/2}_2 
= p_1^2 + 2p_2p_1q_1 - q_1^5 + p_2^2q_2 + 3q_1^3q_2 - 2q_1q_2^2 \\
  \qquad \qquad + z_1 (q_1^2-q_2) + z_2(z_2q_1 + q_1q_2 - p_2^2),   \\
\end{array} \right.
\label{1-3}
\end{equation}
\begin{equation}
(\text{P}_\text{II-1})_2 \left\{ \begin{array}{l}
\displaystyle H^{7/2+1}_1
 = 2p_1p_2 - p_2^3-p_1q_1^2 + q_2^2 - z_1p_2 + z_2p_1 + 2 \alpha q_1,   \\
\displaystyle H^{7/2+1}_2 
= -p_1^2 + p_1p_2^2 + p_1p_2q_1^2 + 2p_1q_1q_2  \\
\qquad + z_1p_1 + z_2(z_2p_1 - p_1q_1^2 + p_1p_2)
 - \alpha (2p_2q_1 + 2q_2 + 2z_2q_1),  \\
\end{array} \right.
\label{1-4}
\end{equation}
\begin{equation}
(\text{P}_\text{II-2})_2 \left\{ \begin{array}{l}
\displaystyle H^{5}_1 
=p_1p_2 - p_1q_1^2 - 2p_1q_2 + p_2q_1q_2 + q_1q_2^2 + q_2z_1 + z_2 (q_1q_2-p_1)+ \alpha q_1,   \\
\displaystyle H^{5}_2 
=p_1^2 - p_1p_2q_1 + p_2^2q_2 - 2p_1q_1q_2 - p_2q_2^2 + q_1^2q_2^2 \\
\qquad \qquad  + z_1( q_1q_2- p_1) - z_2(p_1q_1+q_2^2 + q_2z_2) + \alpha p_2,   \\
\end{array} \right.
\label{1-5}
\end{equation}
\begin{equation}
(\text{P}_\text{IV})_2 \left\{ \begin{array}{l}
\displaystyle H^{4+1}_1 
=p_1^2  + p_1p_2 - p_1q_1^2 + p_2q_1q_2- p_2q_2^2 - z_1p_1 + z_2p_2q_2  + \alpha q_2 + \beta q_1,  \\
\displaystyle H^{4+1}_2 
= p_1p_2q_1 - 2p_1p_2q_2 - p_2^2 q_2 + p_2q_1q_2^2 \\
\qquad   + p_2q_2 z_1+ z_2(p_1p_2 - p_2q_2^2 + p_2q_2z_2) + (p_1-q_1q_2 + q_2z_2) \alpha - \beta p_2 ,   \\
\end{array} \right.
\label{1-6}
\end{equation}
respectively, with arbitrary parameters $\alpha, \beta \in \C$
(these notations for Hamiltonian functions are related to the spectral type of 
a monodromy preserving deformation \cite{KNS}).
The weighted degrees of these hierarchies determined by the Newton diagrams are shown in Table 2
(see also Table 3).
From Table 1, 2 and the equations, we deduce the following properties.
\begin{itemize}
\item $\mathrm{deg} (q_i) + \mathrm{deg} (p_i) = \mathrm{deg}(H_1)-1$.
\item $\mathrm{deg}(z_1) = \mathrm{deg}(H_1)-2 $.
\item $\mathrm{deg}(z_i)+\mathrm{deg}(H_i)$ is independent of $i=1,\cdots ,n$.
\item $\displaystyle \min_{1\leq i\leq n}\{ \mathrm{deg} (q_i), \mathrm{deg}(p_i) \} = 1\,\, \text{or}\,\, 2$.
\item The equation (\ref{1-2}) is invariant under the $\Z_s$-action
\begin{eqnarray*}
(q_i, p_i, z_i) \mapsto (\omega ^{\mathrm{deg}(q_i)} q_i,\,\, \omega ^{\mathrm{deg}(p_i)} p_i
,\,\, \omega ^{\mathrm{deg}(z_i)} z_i),
\end{eqnarray*}
where $s:= \mathrm{deg}(H_1) - 1$ and $\omega :=e^{2\pi i/s}$.
\item The symplectic form $\displaystyle \sum^n_{i=1} dq_i \wedge dp_i + \sum^n_{i=1} dz_i \wedge dH_i$
is also invariant under the same $\Z_s$-action, for which $H_i \mapsto \omega ^{\mathrm{deg}(H_i)} H_i$.
\end{itemize}

We decompose the Hamiltonian function $H_i$ into the principle part $H_i^P$
and the non-principle part $H_i^N$.
Then, we further deduce
 \begin{itemize}
\item $H_i^N$ consists of monomials including arbitrary parameters.
\item $\mathrm{deg} (H_i^N) = \mathrm{deg}(H_i) - \mathrm{deg}(H_1) + 1$.
In particular $\mathrm{deg} (H_1^N) =1$.
\item The variety defined by
\begin{eqnarray*}
H_1^P (q_1, \cdots ,q_n, p_1, \cdots ,p_n, 0, \cdots ,0) = 0
\end{eqnarray*}
in $\C^{2n}$ has a unique singularity at the origin.
\end{itemize}

In section 2.4, several above properties will be proved from the others.
For  $(\text{P}_\text{I})$, $(\text{P}_\text{II})$ and $(\text{P}_\text{IV})$,
we have
\begin{eqnarray*}
H_\text{I}^P(q,p,0) = \frac{1}{2}p^2 - 2q^3, \\
H_\text{II}^P(q,p,0) =\frac{1}{2}p^2 - \frac{1}{2}q^4,  \\
H_\text{IV}^P(q,p,0) =  -pq^2 + p^2q.
\end{eqnarray*}
They define $A_2, A_3$ and $D_4$ singularities at the origin, respectively.
In singularity theory, it is known that if a singularity defined by a quasihomogeneous
polynomial $H(x_1, \cdots ,x_n) = 0$ is isolated, then the rational function
\begin{equation}
\chi (T) := \frac{(T^{h-a_1}-1)\cdots (T^{h-a_n}-1)}
{(T^{a_1}-1)\cdots (T^{a_n}-1) }
\label{1-8}
\end{equation}
becomes a polynomial (Poincar\'{e} polynomial), where 
$\mathrm{deg}(x_i) = a_i$ and $\mathrm{deg}(H) = h$.

Motivated by these observation, we classify regular weights $(a_1,\cdots ,a_n, b_1,\cdots ,b_n; h)$
satisfying certain conditions in Section 2.
In particular, for $n=1$ and $2$, we will show that there is a corresponding Painlev\'{e} 
equation for each weight such that 
$\mathrm{deg}(q_i) = a_i, \mathrm{deg}(p_i) = b_i$ and $\mathrm{deg}(H) = h$.
In Section 2.4, a Hamiltonian system, whose Hamiltonian function
satisfies certain assumptions on the quasihomogeneity, will be considered.
Then, some of the above properties of weights will be proved.

In Section 3, a brief review of Kovalevskaya exponents of quasihomogeneous vector fields,
which seems to be closely related to regular weights, is given. 
A list of Kovalevskaya exponents of $4$-dim Painlev\'{e} equations are shown in Table 4.
From the table, it is expected that Painlev\'{e} equations defined by semi-quasihomogeneous  
Hamiltonian functions can be classified by their weights and Kovalevskaya exponents.

In Section 4, Kovalevskaya exponents of quasihomogeneous systems are further studied 
by means of singularity theory and dynamical systems theory.
In general, the level surface of quasihomogeneous Hamiltonian functions
has a singularity at the origin.
The weighted blow-up of the singularity at the origin induces a vector field on the exceptional divisor.
Then, Laurent series solutions, Kovalevskaya exponents and the level surface are investigated
via the vector field.
In particular, it is shown that there is a one-to-one correspondence between 
Laurent series solutions and fixed points of the vector field,
and the eigenvalues of the Jacobi matrix of the vector field at the fixed point precisely 
coincide with Kovalevskaya exponents.
With the aid of these results, it is shown for several $4$-dim Painlev\'{e} equations
that the level surface of Hamiltonian functions can be decomposed into a 
disjoint sum of stable manifolds of the fixed points.

\begin{table}[h]
\begin{center}
\begin{tabular}{|c||c|c|c|}
\hline
 & $\mathrm{deg}(q,p,z)$ & $\mathrm{deg} (H)$ & $\kappa$ \\ \hline \hline
$\text{P}_\text{I}$  & $(2,3,4)$ & 6 & 6  \\ \hline
$\text{P}_\text{II}$  & $(1,2,2)$ & 4 & 4  \\ \hline
$\text{P}_\text{IV}$  & $(1,1,1)$ & 3 & 3  \\ \hline
$\text{P}_{\text{III} (D_8)}$  & $(-1,2,4)$ & 2 & 2\\ \hline
$\text{P}_{\text{III} (D_7)}$  & $(-1,2,3)$ & 2 & 2 \\ \hline
$\text{P}_{\text{III} (D_6)}$  & $(0,1,2)$ & 2 & 2  \\ \hline
$\text{P}_\text{V}$  & $(1,0,1)$ & 2 & 2  \\ \hline
$\text{P}_\text{VI}$  & $(1,0,0)$ & 2 & 2  \\ \hline
\end{tabular}
\end{center}
\caption{ $\mathrm{deg} (H)$ denotes the weighted degree of the Hamiltonian function with
respect to the weight $\mathrm{deg}(q,p,z)$.
$\kappa$ denotes the Kovalevskaya exponent defined in Section 3.}
\end{table}

\begin{table}[h]
\begin{center}
\begin{tabular}{|c||c|c|c|}
\hline
 & $\mathrm{deg}(q_j, p_j)$ & $\mathrm{deg} (z_i)$ & $\mathrm{deg} (H_i)$ \\ \hline \hline
$(\text{P}_\text{I})_n$  & $(2j, 2n+3-2j)$ & $2n-2i+4$ & $2n+2i+2$  \\ \hline
$(\text{P}_\text{II-1})_n$  & $(2j-1, 2n+2-2j)$ & $2n-2i+2$ & $2n+2i$  \\ \hline
$(\text{P}_\text{II-2})_n$  & $(j,n+2-j)$ & $n-i+2$ & $n+i+2$  \\ \hline
$(\text{P}_\text{IV})_n$  & $(j,n+1-j)$ & $n-i+1$ & $n+i+1$ \\ \hline
\end{tabular}
\end{center}
\caption{ Weights for four classes of the Painlev\'{e} hierarchies.}
\end{table}

\begin{table}[h]
\begin{center}
\begin{tabular}{|c||c|c|c|}
\hline
 & $\{ \mathrm{deg}(q_j), \mathrm{deg}(p_j)\}$ & $\mathrm{deg} (z_i)$ & $\mathrm{deg} (H_i)$ \\ \hline \hline
$(\text{P}_\text{I})_2$  & $(2,3,4,5)$ & $6,4$ & $8,10$  \\ \hline
$(\text{P}_\text{I})_3$  & $(2,3,4,5,6,7)$ & $8,6,4$ & $10,12,14$  \\ \hline \hline
$(\text{P}_\text{II-1})_2$  & $(1,2, 3,4)$ & $4,2$ & $6,8$  \\ \hline
$(\text{P}_\text{II-1})_3$  & $(1,2, 3,4,5,6)$ & $6,4,2$ & $8,10,12$  \\ \hline\hline
$(\text{P}_\text{II-2})_2$  & $(1,2,2,3)$ & $3,2$ & $5,6$  \\ \hline
$(\text{P}_\text{II-2})_3$  & $(1,2,2,3,3,4)$ & $4,3,2$ & $6,7,8$  \\ \hline\hline
$(\text{P}_\text{IV})_2$  & $(1,1,2,2)$ & $2,1$ & $4,5$ \\ \hline
$(\text{P}_\text{IV})_3$  & $(1,1,2,2,3,3)$ & $3,2,1$ & $5,6,7$ \\ \hline
\end{tabular}
\end{center}
\caption{ Weights for four classes of the Painlev\'{e} hierarchies when $n=2, 3$,
where $\mathrm{deg}(q_j), \mathrm{deg}(p_j)$'s are shown in ascending order.}
\end{table}

%%%%%%%%%%%%%%%%%%%%%%%%%%%%%%%%%%%%%%%%%%%%%%%%%%%%%%%%%%%%%%%%%%%%%%%%

\section{Classification of regular weights}

Let $a_1 ,\cdots ,a_n, b_1,\cdots ,b_n$ and $h$ be positive integers such that $1 \leq a_i, b_i < h$.
Motivated by the observation in Section 1, we suppose the following.
\\[0.2cm]
\textbf{(W1)} $\displaystyle \min_{1\leq i \leq n} \{ a_i, b_i\} = 1 \,\, \text{or} \,\, 2$.
\\
\textbf{(W2)} $a_i + b_i = h-1$ for $i=1,\cdots ,n$.
\\
\textbf{(W3)} A function 
\begin{equation}
\chi (T) = \frac{(T^{h-a_1}-1)(T^{h-b_1}-1)\cdots (T^{h-a_n}-1)(T^{h-b_n}-1) }
{(T^{a_1}-1)(T^{b_1}-1)\cdots (T^{a_n}-1) (T^{b_n}-1)}
\label{2-1}
\end{equation}
is polynomial.

In Saito\cite{Sai}, a tuple of integers $(a_1,\cdots, a_n, b_1,\cdots  ,b_n; h)$ satisfying (W3) is called a regular weight.
In this paper, a tuple is called a regular weight if it satisfies (W1) to (W3).
In this section, we will classify all regular weights for $n=1,2,3$.
In particular, for $n=1$ and $n=2$, we will show that there are Hamiltonians of Painlev\'{e} equations
associated with regular weights such that 
$\mathrm{deg}(q_i) = a_i, \mathrm{deg}(p_i) = b_i$ and $\mathrm{deg}(H) = h$.

%%%%%%%%%%%%%%%%%%%%%%%%%%%%%%%%%%%%%%%%%%%%%%%%%%%%%%%%%%%%%%%%%%%%%%%%%%%%%%%%%%%%%%%%%%%%%%%%%%

\subsection{$n=1$}

\noindent \textbf{Proposition \thedef.}
When $n=1$, regular weights satisfying (W1) to (W3) are only
\begin{eqnarray*}
(a, b ; h) = (2,3; 6), \quad (1,2;4), \quad (1,1;3).
\end{eqnarray*}
They coincide with the weights $(\mathrm{deg}(q), \mathrm{deg}(p); \mathrm{deg}(H))$ 
of $H_{\text{I}}, H_{\text{II}}$ and $H_{\text{IV}}$ for $2$-dim Painlev\'{e} equations, respectively, given in Sec.1.
\\

Hence, there is a one to one correspondence between regular weights and the $2$-dim Painlev\'{e} 
equations written in polynomial Hamiltonians.
Note that $\mathrm{deg}(z)$ is recovered by the rule $\mathrm{deg}(z) = \mathrm{deg}(H)-2$
(see also Prop.2.5).
Now we show that $H_{\text{I}}, H_{\text{II}}$ and $H_{\text{IV}}$ can be reconstructed from 
the regular weights with the aid of singularity theory.
\\[0.2cm]
\textbf{Step 1.}
Consider generic polynomials $H(q,p)$ whose weighted degrees are
$\mathrm{deg}(q,p;H) = (2,3;6), (1,2;4)$ and $(1,1;3)$.
They are given by
\begin{eqnarray*}
H &=& c_1 p^2 + c_2 q^3, \\
H &=& c_1 p^2 + c_2q^2p + c_3q^4, \\
H &=& c_1 q^3 + c_2 pq^2 + c_3 p^2q + c_4p^3,
\end{eqnarray*}
respectively, with arbitrary constants $c_1, \cdots ,c_4$.
\\
\textbf{Step 2.} Simplify by symplectic transformations.
One of the results are
\begin{eqnarray*}
H &=& \frac{1}{2}p^2 - 2 q^3, \\
H &=& \frac{1}{2}p^2 - \frac{1}{2}q^4, \\
H &=& -pq^2 + p^2 q,
\end{eqnarray*}
respectively.
\\
\textbf{Step 3.} Consider the versal deformations of them\cite{Arn}.
We obtain
\begin{eqnarray*}
H &=& \frac{1}{2}p^2 - 2 q^3 + \alpha _4 q + \alpha _6, \\
H &=& \frac{1}{2}p^2 - \frac{1}{2}q^4 + \alpha _2 q^2 + \alpha _3 q + \alpha _4, \\
H &=& -pq^2 + p^2 q + \alpha _1 pq + \alpha _2 p + \beta_2 q + \alpha _3,
\end{eqnarray*}
respectively, where $\alpha _i, \beta_i \in \C$ are deformation parameters.
The subscripts $i$ of $\alpha _i, \beta_i$ denote the weighted degrees of $\alpha _i, \beta_i$
so that $H$ becomes a quasihomogeneous.
\\
\textbf{Step 4.} Now we use the ansatz $\mathrm{deg}(z) = \mathrm{deg}(H)-2$ observed in Sec.1.
If there is a parameter $\alpha _i$ such that $i= \mathrm{deg}(H)-2$, then replace it by $z$.
The results are 
\begin{eqnarray*}
H &=& \frac{1}{2}p^2 - 2 q^3 + z q + \alpha _6, \\
H &=& \frac{1}{2}p^2 - \frac{1}{2}q^4 + z q^2 + \alpha _3 q + \alpha _4, \\
H &=& -pq^2 + p^2 q + z pq + \alpha _2 p + \beta_2 q + \alpha _3,
\end{eqnarray*}
respectively.
They are equivalent to $H_{\text{I}}, H_{\text{II}}$ and $H_{\text{IV}}$ up to the scaling of $z$
(constant terms in Hamiltonians such as $\alpha _6$ do not play a role).

Hence, when $n=1$, there is a one to one correspondence between the regular weights and 2-dim
polynomial Painlev\'{e} equations, and we can recover one of them from the other.

%%%%%%%%%%%%%%%%%%%%%%%%%%%%%%%%%%%%%%%%%%%%%%%%%%%%%%%%%%%%%%%%%%%%%%%%%%%%%%%%%%%%%%%%%%%%%%%%%%

\subsection{$n=2$}

\noindent \textbf{Proposition \thedef.}
When $n=2$, regular weights satisfying (W1) to (W3) are only
\begin{eqnarray*}
(a_1, a_2, b_2, b_1 ; h) &=& (2,3,4,5; 8), \\
 &=& (1,2,3,4; 6), \\
 &=& (2,2,3,3; 6), \\
 &=& (1,2,2,3; 5), \\
 &=& (1,1,2,2; 4), \\
 &=& (1,1,1,1; 3),
\end{eqnarray*}
where we assume without loss of generality that $a_1 \leq a_2 \leq b_2 \leq b_1$.
For each weight, there exists a polynomial Hamiltonian of a $4$-dim Painlev\'{e} equation (not unique).
Explicit forms of Hamiltonian functions are given as follows.
\\[0.2cm]
\textbf{(2,3,4,5;8).} The first Hamiltonian $H^{9/2}_1$ of $(\text{P}_\text{I})_2$ shown in Eq.(\ref{1-3})
has this weight with $\mathrm{deg}(q_1, q_2, p_1, p_2) = (2,4,5,3)$.
Another example is
\begin{equation}
H_{\text{Cosgrove}} = -4p_1p_2 - 2p_2^2q_1 - \frac{73}{128}q_1^4 + \frac{11}{8}q_1^2q_2
 - \frac{1}{2}q_2^2- q_1z - \frac{1}{48}\left( q_1 + \frac{\alpha }{6} \right) q_1^2\alpha .
\label{2-2}
\end{equation}
This Hamiltonian system is derived by a Lie-algebraic method of type $B_2$ and 
can be written in Lax form \cite{Chi4}, thus it enjoys the Painlev\'{e} property.
It seems that it does not appear in the list of $4$-dim Painlev\'{e} equations in \cite{K, KNS, Na}.
If we rewrite the system as the fourth order single equation of $q_1 = y$, we obtain
\begin{equation}
y'''' = 18 yy'' + 9(y')^2 - 24y^3 + 16z + \alpha y (y + \frac{1}{9}\alpha ).
\end{equation}
This equation was first given in Cosgrove \cite{Cos}, denoted by F-VI.
He conjectured that this equation defines a new Painlev\'{e} transcendents
(i.e. it is not reduced to known equations).
\\[0.2cm]
\textbf{(1,2,3,4;6).} The first Hamiltonian $H^{7/2+1}_1$ of $(\text{P}_\text{II-1})_2$ shown in Eq.(\ref{1-4})
has this weight $\mathrm{deg}(q_1, q_2, p_1, p_2) = (1,3,4,2)$.
Another example is the matrix Painlev\'{e} equation of the first type $H^{\mathrm{Mat}}_\text{I}$\cite{K, KNS}
defined by
\begin{equation}
H^{\mathrm{Mat}}_\text{I} = \frac{1}{2}p_1^2 - 2q_1^3 - 2p_2^2q_2 + 6q_1q_2 - 2q_1z + 2\alpha p_2,
\label{2-4}
\end{equation}
with $\mathrm{deg}(q_1, q_2, p_1, p_2) = (2,4,3,1)$.
\\[0.2cm]
\textbf{(2,2,3,3;6).} For $H^{7/2+1}_1$ and $H^{7/2+1}_2$ of $(\text{P}_\text{II-1})_2$ shown in Eq.(\ref{1-4}), 
perform the symplectic transformation
\begin{equation}
q_1 = - \frac{y_1}{2x_1},\,\, p_1 = -x_1^2,\,\, q_2 = \frac{y_2}{2},\,\, p_2 = 2 x_2.
\end{equation}
Then we obtain the Hamiltonians
\begin{equation}
\left\{ \begin{array}{l}
\displaystyle H_1^{(2,3,2,3)} = -4x_1^2 x_2 - 8x_2^3 + \frac{y_1^2}{4} + \frac{y_2^2}{4} 
        -2 z_1 x_2 - z_2 x_1^2 - \frac{\alpha y_1}{x_1},  \\
\displaystyle H_2^{(2,3,2,3)} = -x_1^4 -4x_1^2x_2^2 - \frac{x_2y_1^2}{2} +\frac{x_1y_1y_2}{2} \\
\displaystyle  \qquad \quad  - z_1x_1^2 - z_2^2x_1^2 - 2z_2x_1^2 x_2 + \frac{z_2y_1^2}{4}
  - \frac{\alpha z_2 y_1}{x_1}+\frac{2\alpha x_2 y_1}{x_1}- \alpha y_2.  \\
\end{array} \right.
\label{2-6}
\end{equation}
Thus, putting $\alpha  =0$ yields semi-quasihomogeneous Hamiltonians of
\\
$\mathrm{deg}(H_1^{(2,3,2,3)}, H_2^{(2,3,2,3)}) = (6,8)$ with respect to
$\mathrm{deg}(x_1, y_1, x_2, y_2) = (2,3,2,3)$ and $\mathrm{deg}(z_1, z_2) = (4,2)$.
Although this is equivalent to $(\text{P}_\text{II-1})_2$ for $\alpha =0$,
they should be distinguished from each other from a view point of a geometric classification
of Painlev\'{e} equations (i.e. a classification based on the spaces of initial conditions)
because the above symplectic transformation is not a one-to-one mapping.
The direct product of two $(\text{P}_\text{I})$ also has this weight, see Example 4.14.
\\[0.2cm]
\textbf{(1,2,2,3;5).} The first Hamiltonian $H^5_1$ of $(\text{P}_\text{II-2})_2$ shown in Eq.(\ref{1-5})
has this weight with $\mathrm{deg}(q_1, q_2, p_1, p_2) = (1,2,3,2)$.
\\[0.2cm]
\textbf{(1,1,2,2;4).} The first Hamiltonian $H^{4+1}_1$ of $(\text{P}_\text{IV})_2$ shown in Eq.(\ref{1-6})
has this weight with $\mathrm{deg}(q_1, q_2, p_1, p_2) = (1,1,2,2)$.
Another example is the matrix Painlev\'{e} equation of the second type $H^{\mathrm{Mat}}_\text{II}$\cite{K, KNS}
defined by
\begin{equation}
H^{\mathrm{Mat}}_\text{II} = \frac{1}{2}p_1^2 - p_1q_1^2 + p_1q_2 - 2p_2^2q_2 - 4p_2q_1q_2 - p_1z
 + 2\alpha p_2 + 2 \beta (p_2 + q_1),
\label{2-7}
\end{equation}
with $\mathrm{deg}(q_1, q_2, p_1, p_2) = (1,2,2,1)$.
The direct product of two $(\text{P}_\text{II})$ also has this weight.
\\[0.2cm]
\textbf{(1,1,1,1;3).} The Noumi-Yamada system of type $A_4$ \cite{KNS, NY} defined by
\begin{equation}
H^{A_4}_{\mathrm{NY}} = 2p_1p_2q_1 + p_1q_1(p_1-q_1-z) + p_2q_2(p_2-q_2-z) + \alpha p_1 + \beta q_1
 + \gamma p_2 + \delta q_2
\label{2-8}
\end{equation}
has the weight $\mathrm{deg}(q_1, q_2, p_1, p_2) = (1,1,1,1)$,
where $\alpha ,\beta, \gamma, \delta $ are arbitrary parameters.
The direct product of two $(\text{P}_\text{IV})$ also has this weight.

%%%%%%%%%%%%%%%%%%%%%%%%%%%%%%%%%%%%%%%%%%%%%%%%%%%%%%%%%%%%%%%%%%%%%%%%%%%%%%%%%%%%%%%%%%%%%%%%%%

\subsection{$n=3$}

To determine all regular weights satisfying (W1) to (W3), the following lemma is useful.
Without loss of generality, we assume $a_1 \leq a_2 \leq \cdots \leq a_n \leq b_n \leq \cdots \leq b_{2} \leq b_1$.
There exist integers $N$ and $j(1), \cdots ,j(N)$ such that 
\begin{eqnarray*}
& & a_1 = \cdots =a_{j(1)} < a_{j(1)+1} = \cdots = a_{j(2)} < \cdots < a_{j(N)+1} = \cdots = a_n \\
& & \quad \leq b_n = \cdots =b_{j(N)+1} < \cdots < b_{j(2)} = \cdots = b_{j(1)+1} < b_{j(1)} = \cdots = b_1.
\end{eqnarray*}
We put $J_l = j(l) - j(l-1) \,\, (l=1,\cdots ,N+1)$, where $j(0) = 0$ and $j(N+1) = n$.
\\[0.2cm]
\noindent \textbf{Lemma \thedef.}
\\
(i) When $N=0$ (i.e. $a_1 = a_n$), then
\begin{eqnarray*}
(a_1, \cdots ,a_n, b_n, \cdots ,b_1; h) &=& (1,\cdots ,1,1,\cdots ,1; 3) \\
&=& (1,\cdots ,1,2,\cdots ,2;4) \\
&=& (2,\cdots ,2,3,\cdots ,3;6).
\end{eqnarray*}
(ii) When $N \geq 1$, the equality $b_{j(i)} = b_{j(i+1)} + 1$ holds for $i=1,\cdots ,N$
and $J_{i+1}\geq J_{i}$ holds for $i=1,\cdots ,N-1$.
If $a_n \neq b_n$, further $b_n = a_n + 1$ and $J_{N+1} \geq J_N$ hold.
\\
(iii) If $a_i < a_{i+1}$ for any $i=1,\cdots ,n-1$, then
\begin{eqnarray*}
(a_1, \cdots ,a_n, b_n, \cdots ,b_1; h) &=& (1,\cdots ,n,n,\cdots ,2n-1; 2n+1) \\
&=& (1,\cdots ,n, n+1,\cdots ,2n; 2n+2) \\
&=& (2,\cdots ,n+1,n+2,\cdots ,2n+1; 2n+4).
\end{eqnarray*}
\\
\textbf{Proof.}
Because of (W2), Eq.(\ref{2-1}) is rewritten as
\begin{equation}
\chi (T) = \frac{(T^{a_1+1}-1) \cdots (T^{a_n+1}-1)(T^{b_n+1}-1)\cdots (T^{b_1+1}-1) }
{(T^{a_1}-1)\cdots (T^{a_n}-1) (T^{b_n}-1) \cdots (T^{b_1}-1)}.
\label{2-9}
\end{equation}

(i) In this case, $a_1 = a_n \leq b_n = b_1$ due to (W2), which implies 
\begin{eqnarray*}
\chi (T) = \frac{(T^{a_1+1} - 1)^n (T^{b_1+1}-1)^n}{(T^{a_1} - 1)^n (T^{b_1}-1)^n}.
\end{eqnarray*}
Since it is polynomial, either $b_1 + 1$ or $a_1+1$ is a multiple of $b_1$.
If $b_1m = b_1 + 1$, then $(m,b_1) = (2,1)$ and we obtain 
$(a_1,\cdots ,a_n, b_n,\cdots ,b_1) = (1,\cdots 1,1,\cdots ,1)$.
If $a_1 = b_1$, the same result is obtained.
Now suppose that $b_1m = a_1 + 1 < b_1 + 1$.
It is easy to verify that $m=1$ and $b_1 = a_1+1$.
Then, 
\begin{eqnarray*}
\chi (T) = \frac{(T^{a_1+2} - 1)^n }{(T^{a_1} - 1)^n}.
\end{eqnarray*}
Since $a_1+2$ is a multiple of $a_1$, we have $a_1m = a_1 + 2$.
This provides $a_1 = 1$ or $2$ (we need not use (W1)).

(ii) In what follows, we suppose that $b_1 >1$.
In this case, $b_j > 1$ for any $j=1,\cdots ,n$ due to the assumption 
$1\leq a_1 \leq \cdots \leq a_n \leq b_n \leq \cdots \leq b_1$ and (W2).
\\[0.2cm]
\textbf{Step 1.} Since $\chi(T)$ is polynomial, there is a multiple of $b_{j(1)}$ among exponents $b_{j(l)}+1$ in the numerator.
If $b_{j(1)}m = b_{j(1)} + 1$, then $(m,b_{j(1)}) = (2,1)$ and it contradicts the assumption $b_{j(1)} = b_1>1$.

If $b_{j(1)}m = b_{j(l)}+1 < b_{j(1)} + 1$ for some $l>1$, it is easy to verify $m=1, l=2$ and $b_{j(1)} = b_{j(2)} +1$.
There are $J_1$ factors $T^{b_{j(1)}}-1$ in the denominator.
This implies that $2J_2 \geq J_1$ when $N=1$ and $a_n = b_n$, and $J_2\geq J_1$ otherwise.
\\[0.2cm]
\textbf{Step 2.} Now we assume that for some $r\leq N$, $b_{j(i)} = b_{j(i+1)} + 1$ holds for $i=1,\cdots ,r-1$.
There exists a multiple of $b_{j(r)}$ among $b_{j(l)} + 1$.
If $l\leq r$, we have
\begin{eqnarray*}
b_{j(r)}m = b_{j(l)} + 1 = b_{j(l+1)} + 2 = \cdots =b_{j(r)} + r-l+1,
\end{eqnarray*}
which yields
\begin{eqnarray*}
1 < b_{j(r)} \leq r-l+1 \leq r.
\end{eqnarray*}
This proves $b_{j(r)} = b_n = a_n = r$ (otherwise, $a_1$ becomes nonpositive).
Hence, $r=N+1$, which contradicts the assumption $r\leq N$.

If $b_{j(r)}m = b_{j(l)} + 1$ for some $l > r$, then $m=1, l=r+1$ and 
$b_{j(r)} = b_{j(r+1)} +1$.
There are $J_r$ factors $T^{b_{j(r)}}-1$ in the denominator.
This implies that $2J_{r+1} \geq J_r$ when $r=N$ and $a_n = b_n$, and $J_{r+1}\geq J_r$ otherwise.
\\[0.2cm]
\textbf{Step 3.} By induction, we obtain $b_{j(i)} = b_{j(i+1)} + 1$ for $i=1,\cdots ,N$,
and $J_{i+1} \geq J_i$ for $i=1,\cdots ,N-1$.
In particular, if $a_n \neq b_n$, $J_{N+1} \geq J_N$ also holds.
\\[0.2cm]
\textbf{Step 4.} 
There exists a multiple of $b_{j(N+1)} = b_n$ among exponents of the numerator.
Suppose $b_{j(N+1)}m = b_{j(l)} + 1$ for some $l = 1,\cdots ,N+1$.
The same argument as Step 2 shows that $a_n = b_n$.
Suppose $b_{j(N+1)}m = a_{j(l)} + 1 < b_{j(N+1)} + 1$ for some $l = 1,\cdots ,N+1$.
Then, we obtain $m=1, l=N+1$ and $b_{j(N+1)} = b_n = a_n + 1$.
This completes the proof of (ii).

(iii) This is verified by a direct calculation with the aid of (ii).  $\Box$
\\[0.2cm]
\noindent \textbf{Proposition \thedef.}
When $n=3$, regular weights satisfying (W1) to (W3) are only
\begin{eqnarray*}
(a_1, a_2, a_3, b_3, b_2, b_1 ; h) &=& (2,3,4,5,6,7; 10), \\
 &=& (2,3,3,4,4,5; 8), \\
 &=& (1,2,3,4,5,6; 8), \\
 &=& (1,2,3,3,4,5; 7), \\
 &=& (2,2,2,3,3,3; 6), \\
 &=& (1,2,2,3,3,4; 6), \\
 &=& (1,1,2,2,3,3; 5), \\
 &=& (1,1,1,2,2,2; 4), \\
 &=& (1,1,1,1,1,1; 3), 
\end{eqnarray*}
where we assume without loss of generality that $a_1 \leq a_2 \leq a_3 \leq b_3\leq b_2 \leq b_1$.

This proposition is easily obtained with the aid of Lemma 2.3.
To find corresponding Painlev\'{e} equations is a future work.
The weights of 6-dim Painlev\'{e} equations $(\text{P}_\text{I})_3, (\text{P}_\text{II-1})_3,
(\text{P}_\text{II-2})_3$ and $(\text{P}_\text{IV})_3$ shown in Table 3 are included in Prop.2.4.
The author does not know a Painlev\'{e} equation
whose Hamiltonian function is semi-quasihomogeneous but its degree does not satisfy (W1) to (W3).

%%%%%%%%%%%%%%%%%%%%%%%%%%%%%%%%%%%%%%%%%%%%%%%%%%%%%%%%%%%%%%%%%%%%%%%%

\subsection{Properties of weights for semi-quasihomogeneous Hamiltonian systems}

We gave the definition of a regular weight which is independent of differential equations so far. 
Now let us consider the $2n$-dimensional Hamiltonian system
\begin{equation}
\frac{dq_i}{dz} = \frac{\partial H}{\partial p_i}, \quad
\frac{dp_i}{dz} = -\frac{\partial H}{\partial q_i}, \quad i=1,\cdots ,n,
\label{3-1c}
\end{equation}
with the Hamiltonian function $H(q_1,\cdots ,q_n, p_1,\cdots ,p_n, z)$.
We suppose the following.
\\[0.2cm]
\textbf{(A1)} $H$ is semi-quasihomogeneous;
it is decomposed into two polynomials as $H = H^P + H^N$.
For the principle part $H^P$,
there exist integers $1\leq a_i, b_i,r < h$ such that
\begin{equation}
H^P(\lambda ^{a}q, \lambda ^bp, \lambda ^rz) = \lambda ^h H^P(q,p,z),
\end{equation}
for any $\lambda \in \C$,
where $\lambda ^aq = (\lambda ^{a_1}q_1 ,\cdots ,\lambda ^{a_n}q_n)$ and 
$\lambda ^bp = (\lambda ^{b_1}p_1 ,\cdots ,\lambda ^{b_n}p_n)$.
\\
\textbf{(A2)} The Hamiltonian vector field of $H^P$ satisfies
\begin{eqnarray*}
\frac{\partial H^P}{\partial p_i}(\lambda ^{a}q, \lambda ^bp, \lambda ^rz)
 = \lambda ^{1+a_i}\frac{\partial H^P}{\partial p_i}(q,p,z),
\quad \frac{\partial H^P}{\partial q_i}(\lambda ^{a}q, \lambda ^bp, \lambda ^rz) 
 = \lambda ^{1+b_i}\frac{\partial H^P}{\partial q_i}(q,p,z).
\end{eqnarray*}
\textbf{(A3)} The non-principle part satisfies
 $H^N(\lambda ^aq,\lambda ^bp,\lambda ^rz) \sim o(\lambda ^h)$ as $|\lambda | \to \infty$.
\\
\textbf{(A4)} The Hamiltonian vector field of $H = H^P+H^N$ is invariant under the $\Z_s$ action
\begin{equation}
(q_j,p_j,z) \mapsto (\omega ^{a_j}q_j, \omega ^{b_j} p_j, \omega ^r z),
\label{3-3c}
\end{equation}
where $s = h-1$ and $\omega := e^{2\pi i/s}$.
\\
\textbf{(A5)} The symplectic form  $\displaystyle \sum^n_{j=1} dq_j \wedge dp_j + dz \wedge dH$
is also invariant under the same $\Z_s$-action, for which $H \mapsto \omega ^{h} H$.
\\

From these assumptions, we will explain some of properties of weights shown in Section 1.
\\[0.2cm]
\textbf{Remark.} 
The assumption (A2) is used to define the Kovalevskaya exponents in the next section.
In this case, we can construct Laurent series solutions of (\ref{3-1c}) systematically.
Due to the assumptions (A1), (A2) and (A5), it is easy to show that the Hamiltonian vector field of $H^P$
is invariant under the action (\ref{3-3c}).
The assumption (A4) requires that the vector field of $H^N$ is also invariant under the action.
Then, Eq.(\ref{3-1c}) induces a rational differential equation on the weighted projective space
$\C P^{2n+1}(a,b,r,s)$ \cite{Chi1, Chi2}.

In what follows, we assume $h\geq 3$ (if $h\leq 2$, Eq.(\ref{3-1c}) is linear).
\\[0.2cm]
 \textbf{Proposition \thedef.}
Suppose that Eq.(\ref{3-1c}) satisfies (A1) to (A5) and $h\geq 3$.
Then,
\\
(i) $a_i + b_i = h-1$ for $i=1,\cdots ,n$,
\\
(ii) $r = h-2$,
\\
(iii) $\mathrm{deg}(H^N) = 1$,
\\
(iv) if Eq.(\ref{3-1c}) is non-autonomous, $\displaystyle \min_{1\leq i\leq n}\{ a_i, b_i\} = \text{1 or 2}.$
\\[0.2cm]
\textbf{Proof.} The first statement (i) immediately follows from (A1) and (A2).

(ii) Because of (A5), there exists an integer $N$ such that $r+h = N(h-1)$.
Since $r<h$, we obtain $0 < r = N(h-1)-h < h$.
This yields $h < N/(N-2)$ if $N\neq 2$.
This contradicts the assumption $h\geq 3$.
Therefore, $N=2$, which proves $r = h-2$.

(iii) Let $q_1^{\mu_1} \cdots q_n^{\mu_n} p_1^{\nu_1}\cdots p_n^{\nu_n}z^\eta$ be a monomial included in $H^N$.
Due to (A3), the exponents satisfy
\begin{eqnarray*}
0\leq \sum^n_{i=1}\left( a_i\mu_i + b_i \nu_i \right) + r\eta \leq h-1.
\end{eqnarray*}
Further, (A4) implies that there exists an integer $N$ such that
\begin{eqnarray*}
\sum^n_{i=1}\left( a_i\mu_i + b_i \nu_i \right) + r\eta - a_j - b_j + r = N(h-1).
\end{eqnarray*}
This and (i),(ii) give
\begin{eqnarray*}
\sum^n_{i=1}\left( a_i\mu_i + b_i \nu_i \right) + r\eta = N(h-1)+1.
\end{eqnarray*}
Hence, we obtain $0\leq N(h-1)+1 \leq h-1$.
This proves $N=0$ and $\sum^n_{i=1}\left( a_i\mu_i + b_i \nu_i \right) + r\eta = 1$. 

(iv) Suppose that $H$ includes $z$.
Since $\mathrm{deg}(H) = h$ and $\mathrm{deg}(z) = h-2$, $z$ is multiplied by
a function whose weighted degree is $2$.
It exists only when $\displaystyle \min_{1\leq i\leq n}\{ a_i, b_i\} = \text{1 or 2}.$ $\Box$

%%%%%%%%%%%%%%%%%%%%%%%%%%%%%%%%%%%%%%%%%%%%%%%%%%%%%%%%%%%%%%%%%%%%%%%%%%%%%%
%%%%%%%%%%%%%%%%%%%%%%%%%%%%%%%%%%%%%%%%%%%%%%%%%%%%%%%%%%%%%%%%%%%%%%%%%%%%%%

\section{Kovalevskaya exponents}

Kovalevskaya exponents are the most important invariants of a quasihomogeneous vector field
related to the Painlev\'{e} test.
Here, we give a brief review of properties of them according to \cite{Chi2}.
Let us consider the system of differential equations on $\C^m$
\begin{equation}
\frac{dx_i}{dz} = f_i(x_1 , \cdots  ,x_m, z) + g_i(x_1 , \cdots  ,x_m,z), \quad i=1, \cdots  ,m,
\label{3-1}
\end{equation}
where $f_i$ and $g_i$ are polynomials in $(x_1, \cdots  ,x_m, z) \in \C^{m+1}$.
We suppose that 
\\[0.2cm]
\textbf{(K1)} $(f_1, \cdots  ,f_m)$ is a quasi-homogeneous vector field satisfying
\begin{equation}
f_i(\lambda ^{a_1}x_1 ,\cdots ,\lambda^{a_m}x_m, \lambda ^r z) = \lambda ^{1+a_i}f_i(x_1 , \cdots ,x_m,z)
\label{3-2}
\end{equation}
for any $\lambda \in \C$ and $i=1, \cdots ,m$, where $(a_1, \cdots  ,a_m, r) \in \Z^{m+1}_{>0}$.
\\[0.2cm]
\textbf{(K2)} $(g_1, \cdots ,g_m)$ satisfies 
\begin{eqnarray*}
g_i(\lambda ^{a_1}x_1 ,\cdots ,\lambda^{a_m}x_m, \lambda ^r z) = o(\lambda ^{a_i+1}), \quad |\lambda | \to \infty.
\end{eqnarray*}

Put $f_i^A (x_1, \cdots ,x_m):= f_i(x_1, \cdots ,x_m, 0)$ and $f_i^{NA} := f_i - f_i^A$
(i.e. $f_i^A$ and $f_i^{NA}$ are autonomous and nonautonomous parts, respectively).
We also consider the truncated system
\begin{equation}
\frac{dx_i}{dz} = f^A_i(x_1 , \cdots  ,x_m), \quad i=1, \cdots  ,m.
\label{3-3}
\end{equation}
By substituting $x_i(z) = c_i(z-z_0)^{-a_i}$ into the truncated system, we find the following definition.
\\[0.2cm]
\textbf{Definition \thedef.}
A root $c = (c_1,\cdots ,c_m)\in \C^m$ of the equation
\begin{equation}
-a_ic_i = f_i^A (c_1, \cdots  ,c_m), \quad i=1, \cdots ,m
\label{3-4}
\end{equation}
is called a balance. 
\\

For each balance,
$x_i(z) = c_i(z-z_0)^{-a_i}$ is an exact solution of the truncated system for any $z_0 \in \C$.
Due to the assumption (K1), $c=0$ is always a balance which corresponds to the fixed point at the origin.
Usually, we assume $c\neq 0$ for a balance.
Considering the variational equation along the exact solution $x_i(z) = c_i(z-z_0)^{-a_i}$ 
suggests the following definition.
\\[0.2cm]
\textbf{Definition \thedef.} For a balance $c = (c_1, \cdots ,c_m) \neq 0$, the matrix 
\begin{equation}
K=K(c):=\Bigl\{ \frac{\partial f_i^A}{\partial x_j}(c_1,\cdots ,c_m)+ a_i\delta _{ij} \Bigr\}^m_{i,j=1}
\label{3-5}
\end{equation}
and its eigenvalues are called the Kovalevskaya matrix and the Kovalevskaya exponents, respectively, of the system (\ref{3-1})
associated with $c$.
\\[0.2cm]
\textbf{Proposition \thedef} (see \cite{Adl, Chi2, Gor2} for the detail.) Suppose (K1) and (K2).
\\
(i)  $-1$ is always a Kovalevskaya exponent with the eigenvector $(a_1c_1, \cdots ,a_mc_m)^T$.
\\
(ii) $\lambda =0$ is a Kovalevskaya exponent associated with $c$ if and only if $c$
is not an isolated root of the equation $-a_ic_i = f_i^A (c_1, \cdots  ,c_m)$.
\\
(iii) The Kovalevskaya exponents are invariant under weight preserving diffeomorphisms.
\\

Consider a formal series solution of Eq.(\ref{3-1}) of the form
\begin{equation}
x_i = c_i(z-z_0)^{-a_i} +  b_{i,1}(z-z_0)^{-a_i+1}+ b_{i,2}(z-z_0)^{-a_i+2} + \cdots 
\label{3-6}
\end{equation}
Coefficients $b_{i,j}$ are determined by substituting it into Eq.(\ref{3-1}).
The column vector $b_j = (b_{1,j} ,\cdots ,b_{m,j})^T$ satisfies
\begin{eqnarray}
(K-jI)b_j = (\text{a function of $c_i$ and $b_{i,k}$ with $k<j$}).
\label{3-7}
\end{eqnarray}
If a positive integer $j$ is not an eigenvalue of $K$, $b_j$ is uniquely determined.
If a positive integer $j$ is an eigenvalue of $K$ and (\ref{3-7}) has no solutions, 
we have to introduce a logarithmic term $\log (z-z_0)$
into the coefficient $b_j$. In this case, the system (\ref{3-1}) has no Laurent series solution
of the form (\ref{3-6}) with a given balance $c$.
If a positive integer $j$ is an eigenvalue of $K$ and (\ref{3-7}) has a solution $b_j$, 
then $b_j + v$ is also a solution for any eigenvectors $v$.
This implies that the series solution (\ref{3-6}) includes a free parameter in $(b_{1,j} ,\cdots ,b_{m,j})$.
Therefore, if (\ref{3-6}) represents a $k$-parameter family of formal Laurent series solutions
which includes $k-1$ free parameters other than $z_0$, at least $k-1$ Kovalevskaya exponents 
have to be nonnegative integers.
Hence, the classical Painlev\'{e} test \cite{Abl, Gor2, Yos2} for the necessary condition 
for the Painlev\'{e} property is stated as follows;
\\[0.2cm]
\textbf{Classical Painlev\'{e} test.} 
If the system (\ref{3-1}) satisfying (K1) and (K2) has the Painlev\'{e} property in the sense that 
any solutions are meromorphic, then there exists a balance $c = (c_1,\cdots ,c_m)$ such that
all Kovalevskaya exponents except for one $-1$ are nonnegative integers
(such a balance is called a principle balance),
and the Kovalevskaya matrix is semisimple.
In this case, (\ref{3-6}) represents an $m$-parameter family of formal Laurent series solutions.
\\

Due to (K1), the system $dx_i/dt = f_i(x_1,\cdots ,x_m, z)$ is invariant under the the $\Z_s$ action
\begin{equation}
(x_1, \cdots ,x_m,z) \mapsto (\omega ^{a_1}x_1, \cdots ,\omega ^{a_m}x_m, \omega ^{r}z), \quad \omega := e^{2\pi i/s},
\label{3-8}
\end{equation}
where $s = r+1$.
We assume that the full system (\ref{3-1}) is also invariant under the same action
(i.e. the perturbation term $g_i$ admit the same $\Z_s$ action as $f_i$);
\\[0.2cm]
\textbf{(K3)} The system (\ref{3-1}) is invariant under the above $\Z_s$ action.
\\[0.2cm]
\textbf{Proposition \thedef.}
Suppose (K1) to (K3).
\textit{If} the system (\ref{3-1}) has a formal series solution (\ref{3-6}),
it is a convergent series on $0 < |z-z_0| < \varepsilon $ for some $\varepsilon >0$.
In particular, when $g_i = 0\,\, (i=1,\cdots ,m)$ this is true without the assumption (K3).
\\

This proposition is shown in \cite{Gor} for autonomous systems and extended to nonautonomous systems
(\ref{3-1}) in \cite{Chi2} by using the weighted projective space $\C P^{m+1}(a_1, \cdots ,a_m, r, s)$.
The assumption (K1) and (K3) are used to confirm that the system (\ref{3-1}) induces a rational 
vector field on $\C P^{m+1}(a_1, \cdots ,a_m, r, s)$.
The classical Painlev\'{e} test gives the \textit{necessary condition} that 
(\ref{3-1}) has an $m$-parameter family of \textit{formal} Laurent series solutions.
Prop.3.4 means that \textit{if} a formal series solution of the form (\ref{3-6}) exists, it is convergent.
In Prop.3.5 of Chiba \cite{Chi2}, the \textit{necessary and sufficient condition} that (\ref{3-1})
has a $k$-parameter family of \textit{convergent} Laurent series solution (\ref{3-6}) is given
under the assumption (K1) to (K3) with the aid of the weighted projective space, 
Kovalevskaya exponents and the normal form theory of dynamical systems.
\\

For the next theorem, we further assume that
\\
\textbf{(S)} A fixed point of the truncated system (\ref{3-3}) is only the origin, i.e,
\begin{equation}
f_i^A(x_1, \cdots ,x_m) = 0 \quad (i=1, \cdots ,m) \Rightarrow (x_1 ,\cdots ,x_m) = (0, \cdots ,0).
\label{3-9}
\end{equation}
\\[0.2cm]
\textbf{Theorem \thedef.\,\cite{Chi2}} 
If the system (\ref{3-1}) satisfies (K1), (K2) and (S), 
any formal Laurent series solutions with a pole at $z=z_0$ are of the form (\ref{3-6}) such that
$(c_1, \cdots ,c_m) \neq (0, \cdots ,0)$.
If we further assume (K3), they are convergent (due to Prop.3.4).
\\

This theorem means that there are no Laurent series solution $(x_1(z) ,\cdots ,x_m(z))$ of (\ref{3-1})
such that the order of a pole of $x_i$ is larger than $a_i$ for some $i$
(For the proof, (S) is essentially used).
Furthermore, if $(c_1 ,\cdots ,c_m) = 0$ 
(i.e. the orders of a pole of $x_1, \cdots ,x_m$ are smaller than $a_1 ,\cdots ,a_m$), 
it should be a local analytic solution.
Therefore, the leading term of a Laurent series solution is strictly given by $c_i (z-z_0)^{-a_i}$
with a given weight $(a_1, \cdots ,a_m)$ and a balance $(c_1, \cdots ,c_m) \neq 0$.
\\

In the rest of this section, we consider the semi-quasihomogeneous Hamiltonian system (\ref{3-1c}).
If it satisfies (A1) to (A5), then it also satisfies (K1) to (K3) and the above results are applicable.
Further, the assumption (S) implies that a singularity of the algebraic variety defined by $\{H=0 \}$ is isolated.
This fact is used to study a relationship between the Painlev\'{e} equations and singularity theory (see Eq.(\ref{1-8})).
The next lemma is well known \cite{Bor, Gor2, HuYa}.
\\[0.2cm]
\textbf{Lemma \thedef.} 
For a semi-quasihomogeneous Hamiltonian system (\ref{3-1c}) of $\mathrm{deg}(H) = h$
satisfying (A1) and (A2),
if $\kappa$ is a Kovalevskaya exponent, so is $\mu$ given by $\kappa +\mu = h-1$.
In particular, $h$ is always a Kovalevskaya exponent for any balances.
\\

\noindent \textbf{Example \thedef.}
The first Painlev\'{e} equation in Hamiltonian form is given by
\begin{eqnarray*}
& & (\text{P}_\text{I}) \left\{ \begin{array}{l}
\displaystyle \frac{dx}{dz} = 6y^2 + z \\[0.2cm]
\displaystyle \frac{dy}{dz} = x,  \\
\end{array} \right.
\end{eqnarray*}
It satisfies the assumptions (A1) to (A5) as is mentioned with $(a_1, a_2 ; h) = (3,2; 6)$ (Table 1).
The balance is uniquely given by $(c_1, c_2) = (-2,1)$.
The associated Laurent series solution is given by
\begin{eqnarray*}
\left(
\begin{array}{@{\,}c@{\,}}
x\\
y
\end{array}
\right) = \left(
\begin{array}{@{\,}c@{\,}}
\! -2 \! \\
0
\end{array}
\right) T^{-3} + 
\left(
\begin{array}{@{\,}c@{\,}}
0 \\
1
\end{array}
\right) T^{-2} - 
\left(
\begin{array}{@{\,}c@{\,}}
\! z_0/5 \!\\
0
\end{array}
\right) T - 
\left(
\begin{array}{@{\,}c@{\,}}
1/2 \\
\! z_0/10 \!
\end{array}
\right) T^2 + 
\left(
\begin{array}{@{\,}c@{\,}}
A_6 \\
-1/6 \!
\end{array}
\right) T^3 + \cdots ,
\end{eqnarray*} 
where $T = z-z_0$ and $A_6$ is an arbitrary constant.
Lemma 3.6 shows that $h=6$ is a Kovalevskaya exponent.
As a result, an arbitrary constant is included in the sixth place from the leading term
(i.e. in the coefficient of $T^{-3 + h} = T^3$).
\\

We give a list of Kovalevskaya exponents of $4$-dim Painlev\'{e} equations shown in Section 2.2.
In Table 4, $H_1^{9/2}, H_1^{7/2+1}, H_1^{5}$ and $H_1^{4+1}$ denote the first Hamiltonians of 
$(\text{P}_\text{I})_2,(\text{P}_\text{II-1})_2, (\text{P}_\text{II-2})_2$ and $(\text{P}_\text{IV})_2$, respectively, given in Section 1
(this notation is related to the spectral type of a monodromy preserving deformation \cite{KNS}).
For example, $(-1,2,3,6) \times 2$ in Table 4 implies that there are two balances $c$, 
for which the associated Kovalevskaya exponents are $\kappa = -1,2,3$ and $6$.
Since Kovalevskaya exponents are invariant under weight preserving diffeomorphisms, 
we can conclude that two Hamiltonian systems having the same weights are actually different systems
if their Kovalevskaya exponents are different from each other.

For a differential equation $dx/dz = f(x, z)$ on $(x, z) \in \C^{m+1}$, 
an $m$-dim manifold $\mathcal{M}(z)$ parameterized by $z$ is called the space of initial conditions 
if any solutions of the system give global holomorphic sections of the 
fiber bundle $\mathcal{P}=\{ (x,z)\, | \, x\in \mathcal{M}(z), z\in \C\}$ over $\C$ \cite{Oka}.
In particular, the space of initial conditions exists for a system having the Painlev\'{e} property 
in the sense that any solutions are meromorphic. 
Many experts believe that the Painlev\'{e} equations can be classified 
by the geometry of the space of initial conditions, which was confirmed for two dimensional 
Painlev\'{e} equations by Sakai \cite{Sak} and Takano et al. \cite{Tak1, Tak2, Tak3}.
In Chiba \cite{Chi2}, an algorithm to construct the space of initial conditions for semi-quasihomogeneous
systems is obtained by the weighted blow-up of the weighted projective space.
The weight for the weighted projective space is just the weight of the variables, and
The weight for the blow-up is given by Kovalevskaya exponents.
This suggests the conjecture that polynomial systems having the Painlev\'{e} property can be classified 
by their weights and Kovalevskaya exponents.

For $2$-dim Painlev\'{e} equations, we have constructed the Painlev\'{e} equations 
$P_{\text{I}}, P_{\text{II}}$ and $P_{\text{IV}}$ from the weights (Prop.2.1).
In this case, the Kovalevskaya exponent is given by $h$ (Lemma 3.6), which is included in
the information of the weight $(a,b; h)$.
For $4$-dim Painlev\'{e} equations listed in Table 4, they are classified by 
the weights with Kovalevskaya exponents.
Thus, the above conjecture looks true at least up to four dimensional quasihomogeneous systems.

As a convenience for readers, we provide a few $4$-dim Painlev\'{e} equations
whose Hamiltonian functions are polynomial but the weights are not positive integers.
Thus, they do not satisfy the assumption (W3).
\begin{eqnarray}
H^{\mathrm{Mat}}_{\text{IV}} &=& 
\frac{p_1^2 q_1}{2} - p_1q_1^2 + p_1q_2 + 2p_1p_2q_2 -4p_2q_1q_2 - 2p_2^2q_1q_2  -p_1q_1z - 2p_2q_2z \nonumber \\
& & \quad + 2p_2q_1\theta_0 - p_1\theta _1 + 2p_2q_1\theta _1-p_1\theta _2 + 2q_1\theta _2 + 2p_2q_1\theta _2, \label{Mat4} \\
H^{(1,2,1,0)} &=& 
 -p_1^2q_1 - 2p_1q_1^2 + 2p_1q_2 - 2p_1p_2q_2 - 2p_2q_1q_2 \nonumber \\
& & \qquad + (2p_1q_1 + 2p_2q_2)z + (2\alpha _2 + 2\beta_2 )q_1 + 2\beta_2 p_1 + 2\beta_3 p_2, \label{1210} \\
H^{(-1,1,4,2)} &=&  p_1 - p_2^2 - 2p_1q_1q_2 - p_2q_2^2 + 2 \beta_3 q_2   + 2\beta_5 q_1 + p_2z. \label{-1142}
\end{eqnarray}
The first one $H^{\mathrm{Mat}}_{\text{IV}}$, whose degree is $\mathrm{deg} (q_1,q_2, p_1,p_2;h) = (1,2,1,0;3)$, 
is the matrix Painlev\'{e} equation of the fourth type $H^{\mathrm{Mat}}_\text{IV}$ \cite{K, KNS}.
$H^{(1,2,1,0)}\,\, (h=3)$ and $H^{(-1,1,4,2)}\,\, (h=4) $ are obtained in \cite{Chi4} by a Lie algebraic method as well as
$H_{\text{Cosgrove}}$.
Although the weight are nonpositive, they still satisfy (W2) and (A1) to (A5).
See also Table 4.

\begin{table}[p]
\begin{center}
\begin{tabular}{|p{1.5cm}||c|p{2.5cm} p{0.6cm}|}
\hline
 & $(a_1,a_2,b_2,b_1; h)$ & $\qquad \kappa$ & \\ \hline \hline
$H_1^{9/2}$ Eq.(\ref{1-3}) & $(2,3,4,5; 8)$ & $(-1,2,5,8)$ $(-3,-1,8,10)$ & $\times 1$ $\times 1$ \\ \hline
$H_{\text{Cosgrove}}$ Eq.(\ref{2-2}) & $(2,3,4,5; 8)$ & $(-1,3,4,8)$ $(-5,-1,8,12)$ & $\times 1$ $\times 1$   \\ \hline
$H_1^{7/2+1}$ Eq.(\ref{1-4}) & $(1,2,3,4; 6)$ & $(-1,2,3,6)$ $(-3,-1,6,8)$  & $\times 2$ $\times 2$ \\ \hline
$H_{\text{I}}^{\text{Mat}}$ Eq.(\ref{2-4}) & $(1,2,3,4; 6)$ & $(-1,2,3,6) $ $(-2,-1,6,7) $ $(-7,-1,6,12) $  & $\times 2$ $\times 1$ $\times 1$ \\ \hline
$H_1^{(2,3,2,3)}$ Eq.(\ref{2-6}) & $(2,2,3,3; 6)$ &  $(-1,1,4,6)$ $(-3,-1,6,8)$  & $\times 1$ $\times 2$ \\ \hline
$H_1^{5}$ Eq.(\ref{1-5}) & $(1,2,2,3; 5)$ & $(-1,1,3,5)$ $(-2,-1,5,6) $  & $\times 2$ $\times 3$  \\ \hline
$H_1^{4+1}$ Eq.(\ref{1-6}) & $(1,1,2,2; 4)$ & $(-1,1,2,4)$ $(-2,-1,4,5)$  & $\times 3$ $\times 5$  \\ \hline
$H_{\text{II}}^{\text{Mat}}$ Eq.(\ref{2-7})  & $(1,1,2,2; 4)$ & $(-1,1,2,4)$ $(-2,-1,4,5)$  
          $(-5,-1,4,8)$ $(-1,-1,4,4)$  & $\times 3$ $\times 2$ $\times 2$ $\times 1$  \\ \hline
$H_{\text{NY}}^{A_4}$ Eq.(\ref{2-8})  & $(1,1,1,1; 3)$ & $(-1,1,1,3)$ $(-1,-1,3,3)$ $(-3,-1,3,5)$  
     & $\times 5$ $\times 5$ $\times 5$ \\ \hline \hline
$H^{\mathrm{Mat}}_{\text{IV}}$ Eq.(\ref{Mat4})  & $(0,1,1,2; 3)$ 
   & $(-1,1,1,3)$ $(-1,-1,3,3)$ $(-2,-1,3,4)$ $(-4,-1,3,6)$
     & $\times 2$ $\times 3$ $\times 2$ $\times 3$ \\ \hline
$H^{(1,2,1,0)}$ Eq.(\ref{1210})  & $(0,1,1,2; 3)$ & $(-1,1,1,3)$ $(-2,-1,3,4)$  
     & $\times 2$ $\times 4$ \\ \hline
$H^{(-1,1,4,2)}$ Eq.(\ref{-1142})  & $(-1,1,2,4; 4)$ & $(-1,1,2,4)$ $(-3,-1,4,6)$ 
     & $\times 2$ $\times 2$ \\ \hline
\end{tabular}
\end{center}
\caption{Weights and Kovalevskaya exponents $\kappa$ of $4$-dim Painlev\'{e} equations.
The weights for dependent variables $q_1,q_2,p_1,p_2$ are shown in ascending order.
For example, $(-1,2,3,6) \times 2$ means that there are two balances whose Kovalevskaya exponents 
are given by $\kappa = -1,2,3,6.$}
\end{table}

% \newpage

%%%%%%%%%%%%%%%%%%%%%%%%%%%%%%%%%%%%%%%%%%%%%%%%%%%%%%%%%%%%%%%%%%%%%%%%%%%%%%%%%%%%%%%%%%%%%%%%%%
%%%%%%%%%%%%%%%%%%%%%%%%%%%%%%%%%%%%%%%%%%%%%%%%%%%%%%%%%%%%%%%%%%%%%%%%%%%%%%%%%%%%%%%%%%%%%%%%%%

\section{Blow-up of quasihomogeneous systems}

Let us investigate the role of Kovalevskaya exponents for quasihomogeneous systems
from a view point of dynamical systems theory.
Since the Kovalevskaya exponents are defined by the autonomous part of a quasihomogeneous system,
we consider the following autonomous system
\begin{equation}
\frac{dx_i}{dz} = f_i(x_1,\cdots ,x_m), \quad i=1,\cdots ,m
\label{4-1}
\end{equation}
satisfying the assumptions (K1) and (S) for the weight $(a_1, \cdots ,a_m) \in \Z^m_{>0}$.
For a balance $c=(c_1, \cdots ,c_m) \neq 0 \in \C^m$ given as a root of $-a_ic_i = f_i(c_1, \cdots ,c_m)$,
$x_i(z) = c_iz^{-a_i}$ is an exact solution.

We introduce the weighted blow-up $\pi : B \to \C^m$ of the system (\ref{4-1}) at the origin 
by the coordinates transformations
\begin{equation}
\left(
\begin{array}{@{\,}c@{\,}}
x_1 \\
x_2 \\
\vdots \\
x_m
\end{array}
\right)
= \left(
\begin{array}{@{\,}l@{\,}}
r_1^{a_1} \\
r_1^{a_2} X_2^{(1)} \\
\vdots \\
r_1^{a_m} X_m^{(1)}
\end{array}
\right) = \left(
\begin{array}{@{\,}l@{\,}}
r_2^{a_1} X_1^{(2)} \\
r_2^{a_2} \\
\vdots \\
r_2^{a_m} X_m^{(2)}
\end{array}
\right)=\cdots = \left(
\begin{array}{@{\,}l@{\,}}
r_m^{a_1} X_1^{(m)} \\
r_m^{a_2} X_2^{(m)} \\
\vdots \\
r_m^{a_m}
\end{array}
\right),
\label{4-2}
\end{equation}
and the blow-up space $B$ by
\begin{eqnarray*}
B = B_1 \cup B_2 \cup \cdots \cup B_m, \quad B_j \simeq \C^m/\Z_{a_j}.
\end{eqnarray*}
Here, the space $\C^m/\Z_{a_j}$ is defined as follows:
Let $(r_1, X_2^{(1)}, \cdots ,X_m^{(1)})$ be the coordinates of $\C^m$.
Then, $B_1$ is defined as a quotient space by the $\Z_{a_1}$ action
\begin{equation}
(r_1, X_2^{(1)}, \cdots ,X_m^{(1)}) \mapsto 
   (e^{2\pi i/a_1}r_1, e^{-2\pi i a_2/a_1}X_2^{(1)}, \cdots , e^{-2\pi i a_m/a_1}X_m^{(1)}),
\label{4-3}
\end{equation}
and similar for $B_2,\cdots ,B_m$.
Let $\pi : B\to \C^m$ be the surjection defined through (\ref{4-2}).
The exceptional divisor
\begin{equation}
D:= \pi^{-1}(\{ 0\}) = \{ r_1 = 0\} \cup \{ r_2 = 0\} \cup \cdots \cup \{ r_m=0\} \subset B
\end{equation}
is isomorphic to the $m-1$ dimensional weighted projective space $\C P^{m-1}(a_1, \cdots ,a_m)$,
and $\pi|_{B\backslash D} : B\backslash D \to \C^m \backslash \{ 0\}$ is a diffeomorphism.
Since $(c_1, \cdots ,c_m) \neq (0, \cdots ,0)$, we assume $c_1 \neq 0$
and denote the first local coordinates $(r_1, X_2^{(1)}, \cdots ,X_m^{(1)})$
on the chart $B_1$ as $(r, X_2, \cdots ,X_m)$ for simplicity.
In this coordinates, Eq.(\ref{4-1}) is written as
\begin{eqnarray*}
\left\{ \begin{array}{l}
\displaystyle \frac{dr}{dz} = \frac{1}{a_1}r^2 f_1(1, X_2, \cdots ,X_m)  \\[0.2cm]
\displaystyle \frac{dX_i}{dz} = rf_i(1, X_2, \cdots ,X_m)-\frac{a_i}{a_1} rX_i f_1(1, X_2, \cdots ,X_m), \quad
i=2,\cdots ,m.  \\
\end{array} \right.
\end{eqnarray*}
A new independent variable $t$ is defined by $d/dz = r \cdot d/dt$, that results in 
\begin{equation}
\left\{ \begin{array}{l}
\displaystyle \frac{dr}{dt} = \frac{1}{a_1}r f_1(1, X_2, \cdots ,X_m)  \\[0.2cm]
\displaystyle \frac{dX_i}{dt} = f_i(1, X_2, \cdots ,X_m)-\frac{a_i}{a_1}X_i f_1(1, X_2, \cdots ,X_m), \quad
i=2,\cdots ,m.  \\
\end{array} \right.
\label{4-5}
\end{equation}
We regard it as a vector field on $B_1$.
The set $\{ (0, X_2, \cdots ,X_m)\} \subset D$ is an invariant manifold.
\\
 
\noindent \textbf{Lemma \thedef.}
(i) For a balance $(c_1 ,\cdots ,c_m)$ of (\ref{4-1}) with $c_1\neq 0$,
\begin{eqnarray}
(r, X_2, \cdots ,X_m) = (0, c_1^{-a_2/a_1}c_2, \cdots , c_1^{-a_m/a_1}c_m)
\label{4-6}
\end{eqnarray}
is a fixed point of the vector field (\ref{4-5}).
Conversely, for any fixed point $(0, X_2, \cdots ,X_m)$ of (\ref{4-5}) on the divisor, 
there exists a balance $(c_1, \cdots , c_m)$ satisfying (\ref{4-6}).
\\
(ii) For a balance $c$, the exact solution $x_i(z) = c_iz^{-a_j},\, (i=1,\cdots ,m)$ 
written on the blow-up space converges to the fixed point (\ref{4-6}) as $z \to \infty$.
\\[0.2cm]
\textbf{Proof.} (i)
A fixed point satisfying $r=0$ is given by a root of the equation
\begin{eqnarray}
a_1 f_i(1, X_2, \cdots ,X_m)- a_i X_i f_1(1, X_2, \cdots ,X_m)=0,\,\, (i=2,\cdots ,m).
\label{4-7}
\end{eqnarray}
If there is a root $(X_2, \cdots ,X_m)$ satisfying $f_1(1, X_2, \cdots ,X_m) = 0$,
then $f_i(1, X_2, \cdots ,X_m) = 0$ for $i=2,\cdots ,m$.
This contradicts with the assumption (S).
Thus, there is a number $\lambda  \neq 0$ such that (\ref{4-7}) is equivalent to
\begin{equation}
\left\{ \begin{array}{l}
f_i(1,X_2,\cdots ,X_m) = -a_iX_i\lambda ^{-1}  \\
f_1(1,X_2,\cdots ,X_m) = -a_1\lambda ^{-1}.  \\
\end{array} \right.
\label{4-8}
\end{equation}
By using the assumption (K1), we rewrite the above equation as
\begin{eqnarray*}
& & \lambda ^{a_i+1} f_i(1, X_2, \cdots ,X_m) \\
& = & f_i(\lambda^{a_1} , \lambda ^{a_2} X_2,\cdots ,\lambda ^{a_m} X_n) \\
&=& -\lambda ^{a_i+1} a_i X_i \lambda ^{-1} = -a_i \lambda ^{a_i}X_i,
\end{eqnarray*}
for $i=2,\cdots ,m$, and 
\begin{eqnarray*}
& & \lambda ^{a_1+1} f_1(1, X_2, \cdots ,X_m) \\
& = & f_1(\lambda^{a_1} , \lambda ^{a_2} X_2,\cdots ,\lambda ^{a_m} X_n) \\
&=& -\lambda ^{a_1+1} a_1 \lambda ^{-1} = -a_1 \lambda ^{a_1}.
\end{eqnarray*}
By putting $\lambda^{a_1} =c_1$ and $\lambda ^{a_i}X_i = c_i$, it turns out that 
(\ref{4-8}) is equivalent to the equation $-a_ic_i = f_i(c_1, \cdots ,c_m)$ to determine a balance.
A proof of (ii) is straightforward. $\Box$
\\

Note that the choice of a branch of $c_1^{-a_j/a_1}$ does not make matters because of the $\Z_{a_1}$ action (\ref{4-3}).
When $c_1 =0$, there are no fixed points in the chart $B_1$ but exists in $B_j$ when $c_j\neq 0$.
In this manner, there is a one-to-one correspondence between balances $c$ and 
fixed points of the vector field induced on the divisor, denoted by $\mathrm{P}(c)$.
If we do not assume (S), there is a fixed point of (\ref{4-5}) on the divisor which results from
not a balance but a fixed point of (\ref{4-1}) other than the origin.
The next proposition associates the Kovalevskaya exponents with the local dynamics around a fixed point
of the vector field.
\\

\noindent \textbf{Proposition \thedef.}
Let $\kappa_1 = -1, \kappa_2 ,\cdots , \kappa_m$ be Kovalevskaya exponents of the system (\ref{4-1})
associated with a balance $c=(c_1, \cdots ,c_m)$.
The eigenvalues of the Jacobi matrix of the vector field (\ref{4-5}) at the fixed point $\mathrm{P}(c)$ are given by
\begin{equation}
\lambda _1 = -c_1^{-1/a_1},\, \lambda _2 = c_1^{-1/a_1} \kappa_2, \, \cdots ,\lambda _m = c_1^{-1/a_1} \kappa_m.
\end{equation}
Hence, the ratio of eigenvalues are the same as that of the Kovalevskaya exponents.
\\[0.2cm]
\textbf{Proof.}
Let $K$ be the Kovalevskaya matrix for a balance $c$.
Set $v_1 = a_1c_1$ and $v_2 = (a_2c_2 ,\cdots ,a_mc_m)$.
Then, $(v_1, v_2)^T$ is an eigenvector of $K$ associated with $\kappa_1 = -1$ (Prop.3.3).
Define
\begin{eqnarray*}
P=\left(
\begin{array}{@{\,}cc@{\,}}
v_1 & 0 \\
v_2 & \mathrm{id}
\end{array}
\right), \quad P^{-1} = \left(
\begin{array}{@{\,}cc@{\,}}
v_1^{-1} & 0 \\
-v_1^{-1} v_2 & \mathrm{id} 
\end{array}
\right), \quad K= \left(
\begin{array}{@{\,}cc@{\,}}
K_1& K_2\\
K_3& K_4
\end{array}
\right).
\end{eqnarray*}
We obtain
\begin{eqnarray*}
P^{-1}KP = \left(
\begin{array}{@{\,}cc@{\,}}
-1 & v_1^{-1}K_2 \\
0 & K_4-v_1^{-1} v_2 K_2
\end{array}
\right) =: \left(
\begin{array}{@{\,}cc@{\,}}
-1 & v_1^{-1}K_2 \\
0 & \tilde{K}
\end{array}
\right),
\end{eqnarray*}
where an $(m-1) \times (m-1)$ matrix $\tilde{K} = (\tilde{K}_{ij})_{i,j=2}^m$ is given by 
\begin{equation}
\tilde{K}_{ij} = \frac{\partial f_i}{\partial x_j}(c) 
  - \frac{a_ic_i}{a_1c_1}\frac{\partial f_1}{\partial x_j}(c) + a_i \delta _{ij}.
\end{equation}
By the definition, eigenvalues of $\tilde{K}$ are $\kappa_2, \cdots ,\kappa_m$.

On the other hand, the Jacobi matrix of (\ref{4-5}) at the fixed point $\mathrm{P}(c)$ is given by
\begin{eqnarray*}
J = \left(
\begin{array}{@{\,}cc@{\,}}
\lambda _1 & 0 \\
0 & \tilde{J}
\end{array}
\right),
\end{eqnarray*}
where
\begin{eqnarray*}
\lambda _1 = \frac{1}{a_1} f_1(1, c_1^{-a_2/a_1}c_2, \cdots , c_1^{-a_m/a_1}c_m)
 = \frac{1}{a_1} c_1^{-(a_1+1)/a_1} f_1(c_1, \cdots ,c_m) = -c_1^{-1/a_1},
\end{eqnarray*}
and 
\begin{eqnarray*}
\tilde{J}_{ij} 
&=& \frac{\partial f_i}{\partial x_j}(1, c_1^{-a_2/a_1}c_2, \cdots , c_1^{-a_m/a_1}c_m)
   - \frac{a_i}{a_1}X_i\frac{\partial f_1}{\partial x_j}(1, c_1^{-a_2/a_1}c_2, \cdots , c_1^{-a_m/a_1}c_m) \\
& &  - \frac{a_i}{a_1}f_1(1, c_1^{-a_2/a_1}c_2, \cdots , c_1^{-a_m/a_1}c_m) \delta _{ij} \\
&=& c_1^{-(a_i+1-a_j)/a_1}\frac{\partial f_i}{\partial x_j}(c) 
   - \frac{a_i}{a_1}c_1^{-a_i/a_1} c_ic_1^{-(a_1+1-a_j)/a_1}\frac{\partial f_1}{\partial x_j}(c)
   - \frac{a_i}{a_1}c_1^{-(a_1+1)/a_1} f_1(c) \delta _{ij} \\
&=& c_1^{-1/a_1} c_1^{-(a_i-a_j)/a_1} \left( \frac{\partial f_i}{\partial x_j}(c)
 - \frac{a_ic_i}{a_1c_1} \frac{\partial f_1}{\partial x_j}(c) + a_i \delta _{ij} \right).
\end{eqnarray*}
This shows
\begin{eqnarray*}
c_1^{1/a_1}\tilde{J}_{ij} = c_1^{-(a_i-a_j)/a_1}\tilde{K}_{ij}.
\end{eqnarray*}
Let $\kappa$ be an eigenvalue of $\tilde{K}$ with the eigenvector $u = (u_2, \cdots ,u_m)^T$
satisfying $\sum \tilde{K}_{ij}u_j = \kappa u_i$.
Putting $u_j = c_1^{a_j/a_1} \tilde{u}_j$ yields
\begin{eqnarray*}
\sum c_1^{-(a_i-a_j)/a_1} \tilde{K}_{ij} \tilde{u}_j = \kappa \tilde{u}_i.
\end{eqnarray*}
This proves that $\kappa$ is an eigenvalue of the matrix $c_1^{1/a_1} \tilde{J}$. $\Box$
\\

We turn to the quasihomogeneous Hamiltonian system of degree $m$
\begin{equation}
\frac{dq_i}{dz} = \frac{\partial H}{\partial p_i} ,
\quad \frac{dp_i}{dz} = -\frac{\partial H}{\partial q_i}, \quad i=1,\cdots ,m.
\label{4-11}
\end{equation}
We assume stronger conditions than (A1) and (A2) as follows.
\\

\noindent \textbf{(H0)} There exist polynomials $H=H_1,H_2 ,\cdots ,H_k\,\, (1\leq k \leq m)$ that commute
with respect to the canonical Poisson structure; $\{ H_i, H_j\} = 0$ for $i,j=1,\cdots ,k$.
\\
\noindent \textbf{(H1)}
$H_i$ is quasihomogeneous; there exist positive integers $a_j, b_j$ and $h_i$ such that
\begin{equation}
H_i(\lambda ^{a}q, \lambda ^bp) = \lambda ^{h_i} H_i(q,p), \quad i=1,\cdots ,k,
\label{4-12}
\end{equation}
for any $\lambda \in \C$, where $\lambda ^aq = (\lambda ^{a_1}q_1 ,\cdots ,\lambda ^{a_m}q_m)$ and 
$\lambda ^bp = (\lambda ^{b_1}p_1 ,\cdots ,\lambda ^{b_m}p_m)$.
\\
\noindent \textbf{(H2)} $h_1 = a_j + b_j + 1$ for $j=1,\cdots ,m$ and $h_1 \leq h_i$ for $i=1,\cdots ,k$.
\\
\noindent \textbf{(S)} A fixed point is only the origin;
\begin{eqnarray*}
\frac{\partial H_1}{\partial p_i} = \frac{\partial H_1}{\partial q_i} = 0 \,\, (i=1,\cdots ,m)
\Rightarrow (q_1, \cdots ,q_m,p_1, \cdots ,p_m) = 0.
\end{eqnarray*}

Note that (H1) and (H2) for $k=1$ is equivalent to (A1) and (A2), see Prop.2.5 (i).
Let $c = (c_1, \cdots ,c_{2m})$ be a balance determined by $H_1$ as
\begin{equation}
\frac{\partial H_1}{\partial p_i}(c) = -a_ic_i, \quad \frac{\partial H_1}{\partial q_i}(c)=b_ic_{m+i},
\quad i=1,\cdots ,m.
\label{4-13}
\end{equation}
The Kovalevskaya matrix at $c$ is 
\begin{eqnarray*}
K(c) =
\left(
\begin{array}{@{\,}cc@{\,}}
\displaystyle \frac{\partial ^2H_1}{\partial p\partial q}(c) & 
  \displaystyle \frac{\partial ^2H_1}{\partial p\partial p}(c) \\[0.2cm]
-\displaystyle \frac{\partial ^2H_1}{\partial q\partial q}(c) & 
 -\displaystyle \frac{\partial ^2H_1}{\partial q\partial p}(c)
\end{array}
\right)
+\left(
\begin{array}{@{\,}cc@{\,}}
\mathrm{diag}(a_1, \cdots ,a_m) & 0 \\[0.5cm]
0 & \mathrm{diag}(b_1, \cdots ,b_m)
\end{array}
\right),
\end{eqnarray*}
where $\partial ^2 /\partial p\partial q
 = \left( \partial ^2 / \partial p_i \partial q_j \right)_{i,j}$.
\\

\noindent \textbf{Lemma \thedef.}
\begin{eqnarray}
\sum^m_{j=1} \left( a_jq_j \frac{\partial H_i}{\partial q_j} + b_jp_j \frac{\partial H_i}{\partial p_j}\right)
 = h_i H_i(q,p).
\label{4-14}
\end{eqnarray}
\textbf{Proof.} This is obtained by the derivative of (\ref{4-12}) at $\lambda =1$.
\\

\noindent \textbf{Lemma \thedef.}
For any balance $c$, $H_i (c) = 0$ for $i=1,\cdots ,k$.
\\[0.2cm]
\textbf{Proof.} Use the relations (\ref{4-13}), (\ref{4-14}) and $\{ H_1, H_j\} = 0$ at $(q,p)=c$.
\\

In what follows, the gradient of a function $H$ is denoted by a row vector
\begin{equation}
dH := \left( \frac{\partial H}{\partial q_1} ,\cdots , \frac{\partial H}{\partial q_m},
 \frac{\partial H}{\partial p_1},\cdots ,\frac{\partial H}{\partial p_m}\right) .
\end{equation}
The assumption (S) implies that $dH_1(q,p) = 0$ if and only if $(q,p)=0$;
i.e. the origin is a unique singularity of the variety $\{H_1 = 0\}$.
The following result was first obtained by Yoshida \cite{Yos1}.
Here, we give a simple proof.
\\

\noindent \textbf{Theorem \thedef.}
For a balance $c$, the following equality
\begin{equation}
dH_j(c) (K(c) - h_j\cdot \mathrm{Id}_{2m\times 2m}) = 0, \quad j=1,\cdots ,k.
\label{4-16}
\end{equation}
holds. In particular, if $dH_j(c) \neq  0$, then $h_j$ is a Kovalevskaya exponent.
\\[0.2cm]
\textbf{Proof.}
$\{ H_1, H_j\} = 0$ gives
\begin{eqnarray*}
0 &=& \frac{\partial }{\partial q_l}\{ H_1, H_j\} \\
&=& \sum^m_{i=1} \left( \frac{\partial ^2 H_1}{\partial q_i \partial q_l} \frac{\partial H_j}{\partial p_i}
    +  \frac{\partial H_1}{\partial q_i}\frac{\partial ^2 H_j}{\partial p_i \partial q_l}
    - \frac{\partial ^2 H_j}{\partial q_i \partial q_l} \frac{\partial H_1}{\partial p_i}
    -  \frac{\partial H_j}{\partial q_i}\frac{\partial ^2 H_1}{\partial p_i \partial q_l} \right).
\end{eqnarray*}
Substituting $(q,p) = c$ yields
\begin{eqnarray*}
0 = \sum^m_{i=1} \left( a_ic_i \frac{\partial ^2 H_j}{\partial q_i \partial q_l} 
    +  b_ic_{m+i}\frac{\partial ^2 H_j}{\partial p_i \partial q_l}
    +  \frac{\partial ^2 H_1}{\partial q_i \partial q_l} \frac{\partial H_j}{\partial p_i}
    -  \frac{\partial H_j}{\partial q_i}\frac{\partial ^2 H_1}{\partial p_i \partial q_l}\right).
\end{eqnarray*}
By the derivative of (\ref{4-14}) with respect to $q_l$, we obtain
\begin{eqnarray*}
\sum^m_{i=1} \left( a_iq_i \frac{\partial ^2 H_j}{\partial q_i\partial q_l}
 + b_ip_i \frac{\partial ^2H_j}{\partial p_i\partial q_l} \right) + a_l \frac{\partial H_j}{\partial q_l}
 = h_j \frac{\partial H_j}{\partial q_l}.
\end{eqnarray*}
Thus, we obtain
\begin{equation}
0 = \sum^m_{i=1}\left( \frac{\partial ^2 H_1}{\partial q_i \partial q_l} \frac{\partial H_j}{\partial p_i}
    -  \frac{\partial ^2 H_1}{\partial p_i \partial q_l}\frac{\partial H_j}{\partial q_i}\right)
 + (h_j - a_l) \frac{\partial H_j}{\partial q_l}
\end{equation}
for $l=1, \cdots ,m$.
The derivative of (\ref{4-14}) with respect to $p_l$ gives similar $m$ relations.
The resultant $2m$ relations are equivalent to (\ref{4-16}). $\Box$
\\

\noindent \textbf{Example \thedef.}
We consider the Hamiltonians of degree 2
\begin{equation}
 \left\{ \begin{array}{l}
\displaystyle H_1 
=2p_2p_1 + 3p_2^2q_1 + q_1^4 - q_1^2q_2 - q_2^2,  \\
\displaystyle H_2 
= p_1^2 + 2p_2p_1q_1 - q_1^5 + p_2^2q_2 + 3q_1^3q_2 - 2q_1q_2^2.
\end{array} \right.
\label{4-18}
\end{equation}
They are autonomous parts of $(\text{P}_\text{I})_2$ given in (\ref{1-3}).
They satisfy (H0), (H1), (H2) and (S) for the weight $(a_1, a_2, b_1, b_2) = (2,4,5,3)$
and $h_1 = 8, h_2 = 10$ shown in Table 2.
There are two balances $c_1=(1,1,1,-1)$ and $c_2 = (3,0,27,-3)$.
For the former $c_1$, the Kovalevskaya exponents are $\kappa = -1,2,5,8$.
Thus, the corresponding Laurent series solution (\ref{3-6}) represents a general solution
including four free parameters (Painlev\'{e} test).
Since $\kappa \neq -10$, Thm.4.5 implies that $dH_2(c_1) =0$.
For the balance $c_2$, we can verify that $dH_1(c_2) \neq 0,\, dH_2(c_2) \neq 0$.
Therefore, Thm.4.5 and Lemma 3.6 show that the Kovalevskaya exponents are given by $\kappa = -3,-1,8,10$
(see $H^{9/2}_1$ of Table 4).
\\

Let $V$ be a variety defined by the level set
\begin{equation}
V = \{ (q,p) \in \C^{2m} \, | \, H_j(q,p) = 0,\, j=1,\cdots ,k\} \ni 0.
\label{4-19}
\end{equation}
Lemma 4.4 shows $c\in V$ for a balance $c$.
Because of (H1), the orbit
\begin{eqnarray*}
\{ (\lambda ^{a_1}c_1 ,\cdots , \lambda ^{a_m}c_m,
 \lambda ^{b_1}c_{m+1},\cdots ,\lambda ^{b_m}c_{2m}) \, | \, \lambda \in \C\}
\end{eqnarray*}
is also included in $V$.
Let us consider the weighted blow-up $\pi : B\to \C^{2m}$ at the origin
\begin{equation}
B = B_1\cup \cdots \cup B_{2m},\quad B_i = \C^{2m}/\Z_{a_i},\, B_{m+i} = \C^{2m}/\Z_{b_i}
\end{equation}
for $i=1,\cdots ,m$.
The exceptional divisor $D$ is a $2m-1$ dimensional weighted projective space
\begin{equation}
D = \pi^{-1}(\{ 0\}) = \C P^{2m-1} (a_1, \cdots ,a_m, b_1,\cdots ,b_m).
\end{equation}
For a balance $c$, we assume $c_1 \neq 0$ as before.
The local coordinates $(r, Q_2, \cdots ,Q_m, P_1, \cdots ,P_m)$ on $B_1$ is defined by
\begin{eqnarray*}
& & q_1 = r^{a_1},\, q_i=r^{a_i}Q_i\quad (i=2,\cdots ,m),\\
& &  p_i = r^{b_i}P_i\quad  (i=1,\cdots ,m).
\end{eqnarray*}
In particular, $D \cap B_1$ is given by the set $\{ r=0\}$.
The set $\pi^{-1} (V) \subset B$ is a disjoint union of $D$ and $\pi^{-1}(V \backslash \{ 0\})$.
Let $\overline{\pi^{-1}(V \backslash \{ 0\})}$ be the closure with respect to the usual topology and 
\begin{eqnarray*}
V_0:= D\cap \overline{\pi^{-1}(V \backslash \{ 0\})} \subset D,
\end{eqnarray*}
see Fig.1.
On the chart $B_1$, we have $H_j(q,p) = r^{h_j} H_j(1, Q_2, \cdots ,P_m)$.
Hence, define
\begin{eqnarray*}
V_{01}:= 
\{ (0, Q_2,\cdots ,P_m)\in D \, | \, H_j(1, Q_2,\cdots ,P_m ) = 0,\,\, j=1,\cdots ,k\} \subset D\cap B_1.
\end{eqnarray*}
The sets $V_{0i}$ on the chart $B_i$ are also defined in the same way for $i=2,\cdots ,2m$.
Then, we have
\begin{eqnarray*}
& & V_0 = V_{01}\cup V_{02} \cup \cdots  \cup V_{0,2m} \subset D \\
& & \pi^{-1}(V) = D\cup \pi^{-1}(V \backslash \{ 0\}) \simeq D\cup (V_0 \times \C).
\end{eqnarray*}
$V$ and $V_0$ are $2m-k$ and $2m-k-1$ dimensional manifolds, respectively, with singularities.
As in Eq.(\ref{4-5}), the system (\ref{4-11}) induces the vector field $\mathcal{X}$ on $B$
after a suitable change of the independent variable.
On $B_1$, $\mathcal{X}$ is expressed as
\begin{equation}
\left\{ \begin{array}{l}
\displaystyle \frac{dr}{dt} = \frac{1}{a_1}r \frac{\partial H_1}{\partial p_1} (1,Q_2, \cdots ,P_m) \\[0.4cm]
\displaystyle \frac{dQ_i}{dt} = \frac{\partial H_1}{\partial p_i}(1,Q_2, \cdots ,P_m)
      - \frac{a_i}{a_1} Q_i \frac{\partial H_1}{\partial p_1}(1,Q_2, \cdots ,P_m)  \\[0.4cm]
\displaystyle \frac{dP_i}{dt} = -\frac{\partial H_1}{\partial q_i}(1,Q_2, \cdots ,P_m)
      - \frac{b_i}{a_1} P_i \frac{\partial H_1}{\partial p_1}(1,Q_2, \cdots ,P_m).
\end{array} \right.
\label{4-22}
\end{equation}

\noindent \textbf{Proposition \thedef.}\\
(i) $D$ is an invariant manifold of $\mathcal{X}$.
\\
(ii) $V_0\subset D$ is an invariant manifold of $\mathcal{X}$.
\\
(iii) All fixed points of $\mathcal{X}$ are included in $V_0$.
\\
(iv) For a balance $c=(c_1, \cdots ,c_{2m})$, the orbit of the exact solution 
$q_i(z) = c_iz^{-a_1}$, $p_i(z) = c_{m+i}z^{-b_i}$ of (\ref{4-11}) is included in $V$.
On the blow-up space $B$, it tends to a fixed point $\mathrm{P}(c)$ on $V_0$ as $z\to \infty$.
\\[0.2cm]
\textbf{Proof.}
It is sufficient to prove the statements on the chart $B_1$.
Since $D\cap B_1 = \{ r=0\}$, (i) immediately follows from Eq.(\ref{4-22}).
By the assumption (S), all fixed points of $\mathcal{X}$ lie on the divisor $D$.
Due to Lemma 4.1(i), a fixed point $\mathrm{P}(c)$ on $D\cap B_1$ is of the form
$(0, c_1^{-a_2/a_1}c_2 ,\cdots ,c_1^{-b_m/a_1}c_{2m})$ for a balance $c$.
Then, Lemma 4.4 implies
\begin{eqnarray*}
H_j(1, c_1^{-a_2/a_1}c_2 ,\cdots ,c_1^{-b_m/a_1}c_{2m})
 = c_1^{-h_j/a_1}H_j(c) = 0,
\end{eqnarray*}
which proves (iii): $\mathrm{P}(c) \in V_{01} \subset V_0$.
Part (iv) is shown by Lemma 4.1 (ii) and Lemma 4.4.
Finally, let us show the statement (ii).
Along an integral curve of (\ref{4-22}), we have
\begin{eqnarray*}
\frac{d}{dt}H_j(1,Q_2,\cdots ,P_m)=
\sum^m_{i=2}\frac{\partial H_j}{\partial q_i} 
    \left( \frac{\partial H_1}{\partial p_i}-\frac{a_i}{a_1}Q_i\frac{\partial H_1}{\partial p_1} \right)
    + \sum^m_{i=1}\frac{\partial H_j}{\partial p_i} 
    \left( -\frac{\partial H_1}{\partial q_i}-\frac{b_i}{a_1}P_i\frac{\partial H_1}{\partial p_1} \right).
\end{eqnarray*}
By introducing a dummy parameter $Q_1=1$, it is rewritten as
\begin{eqnarray*}
\frac{d}{dt}H_j(1,Q_2,\cdots ,P_m)
&=& \sum^m_{i=1}\frac{\partial H_j}{\partial q_i} 
    \left( \frac{\partial H_1}{\partial p_i}-\frac{a_i}{a_1}Q_i\frac{\partial H_1}{\partial p_1} \right)
    + \sum^m_{i=1}\frac{\partial H_j}{\partial p_i} 
    \left( -\frac{\partial H_1}{\partial q_i}-\frac{b_i}{a_1}P_i\frac{\partial H_1}{\partial p_1} \right) \\
&=& \sum^m_{i=1} \left( \frac{\partial H_j}{\partial q_i}\frac{\partial H_1}{\partial p_i}
      - \frac{\partial H_1}{\partial q_i}\frac{\partial H_j}{\partial p_i} \right)
    - \frac{1}{a_1}\frac{\partial H_1}{\partial p_1}
   \sum^m_{i=1}\left( a_iQ_i\frac{\partial H_j}{\partial q_i}+b_iP_i \frac{\partial H_j}{\partial p_i}\right) .
\end{eqnarray*}
Lemma 4.3 and $\{ H_1, H_j\} = 0$ show
\begin{eqnarray*}
\frac{d}{dt}H_j(1,Q_2,\cdots ,P_m)=-\frac{1}{a_1}\frac{\partial H_1}{\partial p_1}h_j H_j(1,Q_2,\cdots ,P_m).
\end{eqnarray*}
This is a linear equation of $H_j(1,Q_2,\cdots ,P_m)$ solved as
\begin{eqnarray*}
H_j(1,Q_2(t),\cdots ,P_m(t)) = H_j(1,Q_2(0),\cdots ,P_m(0)) \cdot 
   \exp \left[ -\frac{h_j}{a_1} \int^t_0 \frac{\partial H_1}{\partial p_1}ds \right].
\end{eqnarray*}
This proves that if $(0,Q_2,\cdots ,P_m) \in V_{01}$ at the initial time $t=0$,
so that \\
$H_j(1,Q_2(0),\cdots ,P_m(0))= 0$,
then $(0,Q_2,\cdots ,P_m) \in V_{01}$ for any $t\in \R$. $\Box$
\\

Fix a balance $c=(c_1, \cdots ,c_{2m}) \neq 0$ with $c_1\neq 1$.
Without loss of generality we assume that $c_1 = 1$ by a suitable scaling of the independent variable.
By Lemma 4.1, the balance associates the fixed point $\mathrm{P}(c):(r, Q_2,\cdots ,P_m) = (0,c_2 ,\cdots ,c_{2m})$
of the vector field (\ref{4-22}) on the chart $B_1$.
Prop.4.2 shows that the Jacobi matrix at the fixed point written in
$(r, Q_2,\cdots ,P_m)$-coordinates is of the form
\begin{equation}
J = \left(
\begin{array}{@{\,}cc@{\,}}
-1 & 0 \\
0 & \tilde{J}
\end{array}
\right),
\end{equation}
and its eigenvalues coincide with the Kovalevskaya exponents $\kappa_1 = -1, \kappa_2,\cdots ,\kappa_{2m}$.
Thus, eigenvectors of $\kappa_2,\cdots ,\kappa_{2m}$ are tangent to the divisor $D=\{ r=0\}$.
Let $\mathbf{E}_s, \mathbf{E}_u$ and $\mathbf{E}_c$ be the stable, unstable and center subspace 
at the point $\mathrm{P}=\mathrm{P}(c)$, which are eigenspaces of eigenvalues with negative real parts,
positive real parts, and zero real parts, respectively.
Let $W_{s}(\mathrm{P}), W_u(\mathrm{P})$ and $W_{c}(\mathrm{P})$ be a local stable manifold,
unstable manifold and center manifold, respectively, 
which are tangent to $\mathbf{E}_s, \mathbf{E}_u$ and $\mathbf{E}_c$ at $\mathrm{P}$.
Because of Lemma 3.6 ($\kappa + \mu = h_1 - 1 >0$), $\mathrm{dim} \mathbf{E}_u \geq m$ 
and $1 \leq \mathrm{dim} \mathbf{E}_s \leq m$.
\\

\noindent \textbf{Proposition \thedef.}
Under the above situation, the unstable manifold $W_u(\mathrm{P})$ is included in $D$
and the stable manifold $W_{s}(\mathrm{P})$ is included in $\overline{\pi^{-1}(V\backslash \{ 0\})}
\simeq V_0 \times \C$.
If there are no purely imaginary eigenvalues ($\neq 0$), 
$W_c(\mathrm{P})$ is included in $V_0$.
\\[0.2cm]
\textbf{Proof}. Since $\mathbf{E}_u$ is tangent to the divisor $D$ and $D$ is an invariant manifold,
$W_u(\mathrm{P}) \subset D$.
Let $x$ be a point on $W_s (\mathrm{P})$ and suppose $x\notin D\cup V$.
Since $x\in W_s(\mathrm{P})$, a solution of (\ref{4-22}) with an initial value at $x$ 
tends to the fixed point $\mathrm{P}$ as $t\to \infty$.
Since $x\notin V$, $H_j (x) \neq 0$ for some $j$.
This is a contradiction because $H_j = 0$ at $\mathrm{P}$ and $H_j$ is a constant along a solution.
Let $x'$ be a point on $W_s (\mathrm{P})$ such that $x'\in D\backslash V_0$.
Then, there is a point $x$ on $W_s (\mathrm{P})$ and $x\notin D \cup V$ because
the eigenvector of $\kappa_1 = -1$ is transverse to $D$, that is again a contradiction. 
If the Kovalevskaya matrix at $c$ has zero eigenvalues, then a balance $c$
is not isolated (Prop.3.3).
Thus, the fixed point $\mathrm{P}(c)$ is not isolated and there exists a neighborhood $U$ 
of $\mathrm{P}(c)$ such that $U\cap V_0$ consists of fixed points of the vector field.
If there are no purely imaginary eigenvalues, $U\cap V_0$ gives a local center manifold.
$\Box$
\\

Next, we consider the system (\ref{4-11}) satisfying (H0), (H1), (H2) and (S) with $k=m$.
In this case, $V$ and $\overline{\pi^{-1}(V\backslash \{ 0\})}$ are $m$ dimensional 
and $V_0$ is an $m-1$ dimensional variety with singularities.
If the system satisfies the Painlev\'{e} property in the sense that any solution
is meromorphic, there is a balance $c$ such that all Kovalevskaya exponents but unique $-1$
are positive integers (Painlev\'{e} test).
Thus, a stable manifold at $\mathrm{P}(c)$ is one dimensional, which is precisely 
given by the orbit of the special solution $q_i(z) = c_i z^{-a_i},\, p_i(z) = c_{m+i} z^{-b_i}$.
The next theorem consider the opposite situation.
\\[0.2cm]
\noindent \textbf{Theorem \thedef.}
Suppose that the system (\ref{4-11}) satisfies (H0), (H1), (H2) and (S) with $k=m$.
Suppose that there exists a balance $c$ such that vectors $dH_1(c), \cdots ,dH_m(c)$ are 
linearly independent.
Then, there exists a neighborhood $U$ of $\mathrm{P}(c)$ such that 
$\overline{\pi^{-1}(V\backslash \{ 0\})} \cup U = W_s(\mathrm{P}(c))\cup U$.
\\[0.2cm]
\textbf{Proof.}
Theorem 4.5 shows that $h_1, \cdots ,h_m >0$ are Kovalevskaya exponents.
Due to Lemma 3.6, negative integers $\mu_j:= h_1 - 1 - h_j,\,\, j=1,\cdots ,m$ are also Kovalevskaya exponents.
Thus, a local stable manifold $W_s(\mathrm{P}(c))$ of the vector field (\ref{4-22}) is an $m$-dimensional
smooth manifold included in $\overline{\pi^{-1}(V\backslash \{ 0\})}$.
Indeed, again Theorem 4.5 implies that the (right) eigenvectors of $K(c)$ associated with eigenvalues
$\mu_1, \cdots ,\mu_m$ are orthogonal to $dH_1(c), \cdots ,dH_m(c)$.
Hence, the stable subspace $\mathbf{E}_s$ coincides with the tangent space of 
$\overline{\pi^{-1}(V\backslash \{ 0\})}$ at $\mathrm{P}(c)$. $\Box$
\\

A balance $c$ satisfying the assumption of the theorem
(i.e. $h_1, \cdots ,h_m >0$ are Kovalevskaya exponents), for which 
$\mathrm{dim} \mathbf{E}_u = \mathrm{dim} \mathbf{E}_s = m$, is called the lowest balance.
The existence of a lowest balance is proved by \cite{Er} for a certain class of 
integrable systems called the hyperelliptically separable systems, while 
the existence for more general systems is not known. 
Let us demonstrate our results for several 4-dim systems obtained from the autonomous parts
of Painlev\'{e} equations. See also Table 5.
They have lowest balances and 
$\overline{\pi^{-1}(V\backslash \{ 0\})}$ is decomposed into the disjoint union
of stable manifolds at the fixed points on the divisor.
\\[0.2cm]
\noindent \textbf{Example \thedef.}
We consider the autonomous part of $(\text{P}_\text{I})_2$ given in Example 4.6.
Since the Kovalevskaya exponents of the balance $c_1=(1,1,1,-1)$, which corresponds to 
the principle Laurent series solution, are $-1,2,5,8$,
the stable manifold $W_s(\mathrm{P}_1)$ at the fixed point $\mathrm{P}_1 = \mathrm{P}(c_1)$ is given by the orbit of the 
special solution $q_1 = z^{-2}, q_2=z^{-4}, p_1=z^{-5}, p_2 = -z^{-3}$.
The Kovalevskaya exponents of the balance $c_2=(3,0,27,-3)$, which satisfies the 
assumptions for Theorem 4.9, are $-3,-1,8,10$.
The $2$-dim stable manifold $W_s(\mathrm{P}_2)$ locally coincides with $\overline{\pi^{-1}(V\backslash \{ 0\})}$.
In this case, $\overline{\pi^{-1}(V\backslash \{ 0\})}$ is decomposed into the disjoint union
of $W_s(\mathrm{P}_1)$ and $W_s(\mathrm{P}_2)$, see Fig.1.

On the $(r, Q_2, P_1, P_2)$-coordinates, $H_1$ and $H_2$ are written as
\begin{eqnarray*}
H_1 &=& r^8 (1+2P_1P_2+3P_2^2-Q_2-Q_2^2), \\
H_2 &=& r^{10} (-1+P_1^2+2P_1P_2 + 3Q_2+P_2^2Q_2-2Q_2^2).
\end{eqnarray*}
Thus, $V_{01}$ is defined by
\begin{eqnarray*}
V_{01}=\{1+2P_1P_2+3P_2^2-Q_2-Q_2^2 = 0,\, -1+P_1^2+2P_1P_2 + 3Q_2+P_2^2Q_2-2Q_2^2 =0\}.
\end{eqnarray*}
Since $dH_2(c_1) = 0$, $V_0$ is singular at $\mathrm{P}_1 : (Q_2, P_1, P_2) = (1,1,-1)$.
We can verify that it is a $A_4$-singularity 
(a singularity whose normal form of defining equation is $y^2 + x^5=0$).
% proved by a suitable transformations.
\\

\begin{figure}
\begin{center}
\includegraphics[scale=0.6]{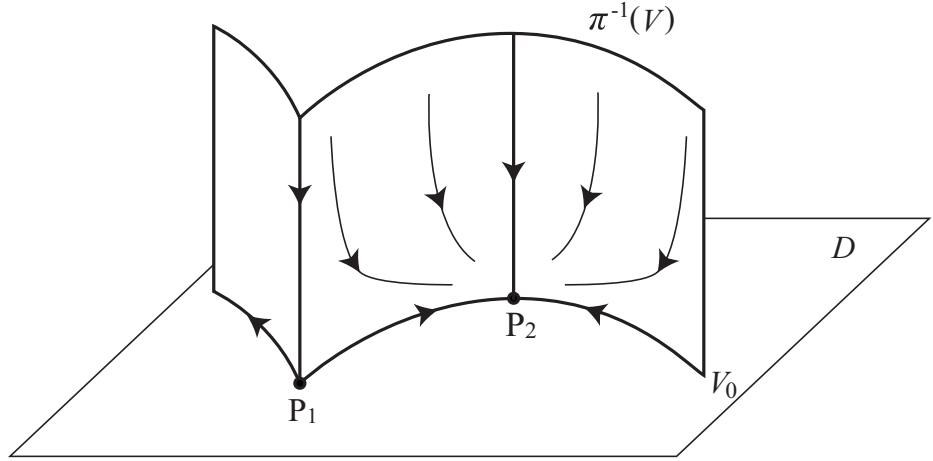}
\caption[]{A schematic view of $\pi^{-1}(V), V_0$ and the dynamics on it for Example 4.10.
The singularity $\mathrm{P}_1$ of $V_0$ is of type $A_4$.}
\label{fig1}
\end{center}
\end{figure}

\noindent \textbf{Example \thedef.}
The autonomous, quasihomogeneous part of $(\text{P}_\text{II-1})_2$ given in (\ref{1-4})
is defined by the Hamiltonians
\begin{equation}
\left\{ \begin{array}{l}
H_1=2p_1p_2 - p_2^3-p_1q_1^2 + q_2^2, \\
H_2=-p_1^2 + p_1p_2^2 + p_1p_2q_1^2 + 2p_1q_1q_2.
\end{array} \right.
\end{equation}
The weights are $(a_1, a_2, b_1, b_2) = (1,3,4,2)$ and $h_1 = 6, h_2 = 8$.
Its four balances and the Kovalevskaya exponents are given by
\begin{eqnarray*}
& & c_1 = (1,0,0,0), \quad \kappa = -1,2,3,6 \\
& & c_2 = (-1,1,0,1), \quad \kappa = -1,2,3,6 \\
& & c_3 = (2,1,0,1), \quad \kappa = -1,-3,6,8 \\
& & c_4 = (-2,3,9,3), \quad \kappa = -1,-3,6,8.
\end{eqnarray*}
Among them, $c_3$ and $c_4$ satisfy the assumptions for Theorem 4.9.
Thus, the fixed points $\mathrm{P}(c_3)$ and $\mathrm{P}(c_4)$ have $2$-dim stable manifolds
that are locally coincide with $\overline{\pi^{-1}(V\backslash \{ 0\})}$.
The fixed points $\mathrm{P}(c_1)$ and $\mathrm{P}(c_2)$ have $1$-dim stable manifolds
that are given by the orbit of special solutions.

On the $(r, Q_2, P_1, P_2)$-coordinates, $H_1$ and $H_2$ are written as
\begin{eqnarray*}
H_1 &=& r^6(2P_1P_2 - P_2^3 - P_1 + Q_2^2) \\
H_2 &=& r^{8} (-P_1^2 + P_1P_2^2 + P_1P_2 + 2P_1Q_2).
\end{eqnarray*}
Thus, $V_{01}$ is defined by
\begin{eqnarray*}
V_{01}=\{2P_1P_2 - P_2^3 - P_1 + Q_2^2 = 0,\, -P_1^2 + P_1P_2^2 + P_1P_2 + 2P_1Q_2 =0\}.
\end{eqnarray*}
Theorem 4.5 shows that $dH_2(c_1) = dH_2(c_2) = 0$.
Hence, $V_0$ is singular at $\mathrm{P}(c_1): (Q_2, P_1, P_2) = (0,0,0)$ and $\mathrm{P}(c_2): (-1,0,1)$.
We can verify that both $\mathrm{P}(c_1)$ and $\mathrm{P}(c_2)$ are $D_5$-singularities
(the normal form is $y(x^2+y^3)=0$).
% $V_0$ is decomposed in to two varieties. $P_1 = 0$ and not.
\\

\noindent \textbf{Example \thedef.}
The autonomous, quasihomogeneous part of $(\text{P}_\text{II-2})_2$ given in (\ref{1-5})
is defined by the Hamiltonians
\begin{equation}
\left\{ \begin{array}{l}
H_1=p_1p_2 - p_1q_1^2 - 2p_1q_2 + p_2q_1q_2 + q_1q_2^2,  \\
H_2=p_1^2 - p_1p_2q_1 + p_2^2q_2 - 2p_1q_1q_2 - p_2q_2^2 + q_1^2q_2^2.
\end{array} \right.
\end{equation}
The weights are $(a_1, a_2, b_1, b_2) = (1,2,3,2)$ and $h_1 = 5, h_2 = 6$.
It has five balances given by
\begin{eqnarray*}
& & c_1 = (1,0,0,0), \quad \kappa = -1,1,3,5 \\
& & c_2 = (-1,-1,1,0), \quad \kappa = -1,1,3,5 \\
& & c_3 = (0,-2,4,-4), \quad \kappa = -1,-2,5,6 \\
& & c_4 = (-2,-2,0,2), \quad \kappa = -1,-2,5,6 \\
& & c_5 = (2,0,0,2), \quad \kappa = -1,-2,5,6.
\end{eqnarray*}
Among them, $c_3, c_4$ and $c_5$ satisfy the assumptions for Theorem 4.9.
Thus, the fixed points $\mathrm{P}(c_3), \mathrm{P}(c_4)$ and $\mathrm{P}(c_5)$ have $2$-dim stable manifolds
that are locally coincide with $\overline{\pi^{-1}(V\backslash \{ 0\})}$.
The fixed points $\mathrm{P}(c_1)$ and $\mathrm{P}(c_2)$ have $1$-dim stable manifolds
that are given by the orbit of special solutions.

On the $(r, Q_2, P_1, P_2)$-coordinates, $H_1$ and $H_2$ are written as
\begin{eqnarray*}
H_1 &=& r^5(P_1P_2 - P_1 - 2P_1Q_2 + P_2Q_2 + Q_2^2) \\
H_2 &=& r^{6} (P_1^2 - P_1P_2 + P_2^2Q_2 - 2P_1Q_2 - P_2Q_2^2 + Q_2^2).
\end{eqnarray*}
$V_{01}$ is defined by
\begin{eqnarray*}
V_{01}&=&\{P_1P_2 - P_1 - 2P_1Q_2 + P_2Q_2 + Q_2^2 = 0, \\
& & \qquad \, P_1^2 - P_1P_2 + P_2^2Q_2 - 2P_1Q_2 - P_2Q_2^2 + Q_2^2 =0\}.
\end{eqnarray*}
Since $dH_2(c_1) = dH_2(c_2) = 0$, $V_0$ is singular at 
$\mathrm{P}(c_1): (Q_2, P_1, P_2) = (0,0,0)$ and $\mathrm{P}(c_2): (-1,-1,0)$.
We can verify that both $\mathrm{P}(c_1)$ and $\mathrm{P}(c_2)$ are $A_5$-singularities
(the normal form is $y^2 + x^6 = 0$).
\\

\noindent \textbf{Example \thedef.}
The autonomous, quasihomogeneous part of $(\text{P}_\text{IV})_2$ given in (\ref{1-6})
is defined by the Hamiltonians
\begin{equation}
\left\{ \begin{array}{l}
H_1=p_1^2  + p_1p_2 - p_1q_1^2 + p_2q_1q_2- p_2q_2^2,  \\
H_2=p_1p_2q_1 - 2p_1p_2q_2 - p_2^2 q_2 + p_2q_1q_2^2.
\end{array} \right.
\end{equation}
The weights are $(a_1, a_2, b_1, b_2) = (1,1,2,2)$ and $h_1 = 4, h_2 = 5$.
It has eight balances given by
\begin{eqnarray*}
& & c_1 = (-1,-1,1,0), \quad \kappa = -1,1,2,4 \\
& & c_2 = (1,0,0,0), \quad \kappa = -1,1,2,4 \\
& & c_3 = (0,1,0,0), \quad \kappa = -1,1,2,4 \\
& & c_4 = (0,-1,2,-4), \quad \kappa = -1,-2,4,5 \\
& & c_5 = (2,0,0,2), \quad \kappa = -1,-2,4,5 \\
& & c_6 = (-1,1,1,0), \quad \kappa = -1,-2,4,5 \\
& & c_7 = (1,2,0,0), \quad \kappa = -1,-2,4,5 \\
& & c_8 = (-2,-2,2,2), \quad \kappa = -1,-2,4,5 
\end{eqnarray*}
Among them, $c_4$ to $c_8$ satisfy the assumptions for Theorem 4.9.
To investigate the fixed points $\mathrm{P}(c_1)$ and $\mathrm{P}(c_2)$, 
we move to $B_1$ chart with the $(r, Q_2, P_1, P_2)$-coordinates, on which $H_1$ and $H_2$ are written as
\begin{eqnarray*}
H_1 &=& r^4(P_1^2  + P_1P_2 - P_1 + P_2Q_2- P_2Q_2^2) \\
H_2 &=& r^{5} (P_1P_2 - 2P_1P_2Q_2 - P_2^2 Q_2 + P_2Q_2^2).
\end{eqnarray*}
$V_{01}$ is defined by
\begin{eqnarray*}
V_{01}&=&\{P_1^2  + P_1P_2 - P_1 + P_2Q_2- P_2Q_2^2 = 0, \\
& & \qquad \, P_1P_2 - 2P_1P_2Q_2 - P_2^2 Q_2 + P_2Q_2^2 =0\}.
\end{eqnarray*}
Since $dH_2(c_1) = dH_2(c_2) = 0$, $V_{01}$ is singular at 
$\mathrm{P}(c_1): (Q_2, P_1, P_2) = (1,1,0)$ and $\mathrm{P}(c_2): (0,0,0)$.
We can verify that both $\mathrm{P}(c_1)$ and $\mathrm{P}(c_2)$ are $D_6$-singularities
(the normal form is $y(x^2+y^4) = 0$).

Note that $\mathrm{P}(c_3)$ is not included in $B_1$ chart because the first component of $c_3$ is zero.
To study $\mathrm{P}(c_3)$ we use to $B_2$ chart with the coordinates $(Q_1, r, P_1, P_2)$
defined like as (\ref{4-2}). In this coordinates, $H_1$ and $H_2$ are written as
\begin{eqnarray*}
H_1 &=& r^4(P_1^2  + P_1P_2 - P_1Q_1^2 + P_2Q_1- P_2) \\
H_2 &=& r^{5} (P_1P_2Q_1 - 2P_1P_2 - P_2^2 + P_2Q_1).
\end{eqnarray*}
$V_{02}$ is defined by
\begin{eqnarray*}
V_{02}&=&\{P_1^2  + P_1P_2 - P_1Q_1^2 + P_2Q_1- P_2 = 0, \\
& & \qquad \, P_1P_2Q_1 - 2P_1P_2 - P_2^2 + P_2Q_1 =0\}.
\end{eqnarray*}
Since $dH_2(c_3)= 0$, $V_{02}$ is singular at $\mathrm{P}(c_3): (Q_1, P_1, P_2) = (0,0,0)$,
which is also a $D_6$-singularity.
\\

\noindent \textbf{Example \thedef.}
Let us consider the following Hamiltonians 
\begin{equation}
\left\{ \begin{array}{l}
H_1= (p_1^2/2 - 2q_1^3) + (p_2^2/2 - 2q_2^3),  \\
H_2=p_1^2/2 - 2q_1^3.
\end{array} \right.
\label{4-27}
\end{equation}
The Hamiltonian equation of $H_1$ is a direct product of the autonomous part of the first Painlev\'{e} equation.
The weights are $(a_1, a_2, b_1, b_2) = (2,2,3,3)$ and $h_1 = 6, h_2 = 6$.
It has three balances $c_1, c_2, c_3$, whose Kovalevskaya exponents are $\kappa = -1,2,3,6$ for $c_1, c_2$ and 
$\kappa = -1,-1,6,6$ for $c_3$.
Since $dH_2(c_1) = dH_2(c_2) = 0$, $V_{0}$ is singular at $\mathrm{P}(c_1)$ and $\mathrm{P}(c_2)$.
We can show that both singularities are $A_2$- singularity.

Similarly, consider the direct product of the autonomous part of the second Painlev\'{e} equation
\begin{equation}
\left\{ \begin{array}{l}
H_1= (p_1^2/2 - q_1^4/2) + (p_2^2/2 - q_2^4/2),  \\
H_2=p_1^2/2 - q_1^4/2.
\end{array} \right.
\label{4-28}
\end{equation}
The weights are $(a_1, a_2, b_1, b_2) = (1,1,2,2)$ and $h_1 = 4, h_2 = 4$.
It has eight balances $c_1,\cdots ,c_8$, whose Kovalevskaya exponents are $\kappa = -1,1,2,4$ for $c_1, c_2, c_3$ and 
$\kappa = -1,-1,4,4$ for the others.
Since $dH_2(c_1) = dH_2(c_2) = dH_2(c_3) = 0$, $V_{0}$ is singular at $\mathrm{P}(c_1)$ to $\mathrm{P}(c_3)$.
We can show that singularities of them are $A_3$- singularity.

Finally, consider the direct product of the autonomous part of the fourth Painlev\'{e} equation
\begin{equation}
\left\{ \begin{array}{l}
H_1= (-p_1q_1^2 + p_1^2q_1) + (-p_2q_2^2 + p_2^2q_2),  \\
H_2=-p_1q_1^2 + p_1^2q_1.
\end{array} \right.
\label{4-29}
\end{equation}
The weights are $(a_1, a_2, b_1, b_2) = (1,1,1,1)$ and $h_1 = 3, h_2 =3$.
It has fifteen balances $c_1,\cdots ,c_{15}$, whose Kovalevskaya exponents are $\kappa = -1,1,1,3$ for $c_1$ to $c_{5}$ and 
$\kappa = -1,-1,3,3$ for the others.
Since $dH_2(c_1) \cdots = dH_2(c_5) = 0$, $V_{0}$ is singular at $\mathrm{P}(c_1)$ to $\mathrm{P}(c_5)$.
They are $D_4$- singularities.

\begin{table}[h]
\begin{center}
\begin{tabular}{|p{1.7cm}||c|p{2.5cm} p{0.6cm}|c|}
\hline
 & $(a_1,a_2,b_1,b_2; h_1, h_2)$ & $\qquad \kappa$ & & \text{sing.} \\ \hline \hline
$H_{\text{I}} \times H_{\text{I}}$ Eq.(\ref{4-27}) & $(2,2,3,3; 6,6)$ & $(-1,2,3,6)$ $(-1,-1,6,6)$ 
  & $\times 2$ $\times 1$ & $A_2$ \\ \hline
$H_{\text{II}} \times H_{\text{II}}$ Eq.(\ref{4-28}) & $(1,1,2,2; 4,4)$ & $(-1,1,2,4)$ $(-1,-1,4,4)$ 
  & $\times 3$ $\times 5$ & $A_3$ \\ \hline
$H_{\text{IV}}\! \times \!H_{\text{IV}}$ Eq.(\ref{4-29}) & $(1,1,1,1; 3,3)$ & $(-1,1,1,3)$ $(-1,-1,3,3)$ 
  & $\times 5$ $\times 10$ & $D_4$ \\ \hline
$H_1^{9/2}$ Eq.(\ref{1-3}) & $(2,4,5,3; 8, 10)$ & $(-1,2,5,8)$ $(-3,-1,8,10)$ & $\times 1$ $\times 1$ & $A_4$ \\ \hline
$H_1^{7/2+1}$ Eq.(\ref{1-4}) & $(1,3,4,2; 6, 8)$ & $(-1,2,3,6)$ $(-3,-1,6,8)$  & $\times 2$ $\times 2$ & $D_5$ \\ \hline
$H_1^{5}$ Eq.(\ref{1-5}) & $(1,2,3,2; 5, 6)$ & $(-1,1,3,5)$ $(-2,-1,5,6) $  & $\times 2$ $\times 3$ & $A_5$  \\ \hline
$H_1^{4+1}$ Eq.(\ref{1-6}) & $(1,1,2,2; 4, 5)$ & $(-1,1,2,4)$ $(-2,-1,4,5)$  & $\times 3$ $\times 5$ & $D_6$  \\ \hline
\end{tabular}
\end{center}
\caption{Weights, Kovalevskaya exponents and singularity-types of $4$-dim autonomous Painlev\'{e} equations.}
\end{table}

%%%%%%%%%%%%%%%%%%%%%%%%%%%%%%%%%%%%%%%%%%%%%%%%%%%%%%%%%%%%%%%%%%%%%%%%%%%%%%%%%%%%%%%%%%%%%%%%%%
%%%%%%%%%%%%%%%%%%%%%%%%%%%%%%%%%%%%%%%%%%%%%%%%%%%%%%%%%%%%%%%%%%%%%%%%%%%%%%%%%%%%%%%%%%%%%%%%%%

\end{document}